\documentclass[11pt]{amsart}
\usepackage[utf8]{inputenc}
\usepackage[english]{babel} 
\usepackage{amsmath} 
\usepackage{amssymb}
\usepackage{amsthm}
\usepackage{mathdots}
\usepackage{mathtools}
\usepackage{enumitem}
\usepackage[vcentermath]{youngtab}
\usepackage{ytableau}
\usepackage{bm} 
\usepackage{eucal} 
\usepackage{hyperref} 
\usepackage{enumitem} 
\usepackage[left=1in,top=1in,right=1in,bottom=1in]{geometry} 
\usepackage{tikz}
\usetikzlibrary{matrix,arrows,decorations.pathmorphing}
\usepackage{fancyhdr} 
\usepackage{parskip} 
\usepackage{microtype} 
\usepackage[scaled=.92]{helvet} 
\usepackage[labelformat=simple,labelfont=rm]{subcaption}

\usepackage{graphicx}
\usepackage{csquotes}
\usepackage{biblatex}
\usepackage{lmodern}
\usepackage[T1]{fontenc}
\usepackage{color}

\title{Parametrizations of $k$-Nonnegative Matrices: Cluster Algebras and $k$-Positivity Tests}
\author{Anna~Brosowsky, Sunita~Chepuri, and Alex~Mason}

\newcommand{\B}[1]{\bm{#1}}

\newcommand{\nth}{^{\text{th}}}

\makeatletter
\def\thm@space@setup{%
  \thm@preskip=\parskip \thm@postskip=0pt
}
\makeatother

\theoremstyle{plain}
\newtheorem{theorem}{Theorem}[section]
\newtheorem*{theorem*}{Theorem}
\newtheorem{conjecture}[theorem]{Conjecture}
\newtheorem{proposition}[theorem]{Proposition}
\newtheorem*{proposition*}{Proposition}
\newtheorem{corollary}[theorem]{Corollary}
\newtheorem*{corollary*}{Corollary}
\newtheorem{lemma}[theorem]{Lemma}
\newtheorem*{lemma*}{Lemma}

\theoremstyle{definition}
\newtheorem*{definition}{Definition}
\newtheorem{example}[theorem]{Example}
\newtheorem{remark}[theorem]{Remark}

\numberwithin{equation}{section}

\newcommand{\Title}{}
\newcommand{\Date}{}
\newcommand{\Name}{}
\newcommand{\Class}{}

\begin{document}
\begin{abstract}
A $k$-positive  matrix is a matrix where all minors of order $k$ or less are positive. 
Computing all such minors to test for $k$-positivity is inefficient, as there are $\sum_{\ell=1}^k \binom{n}{\ell}^2$ of them in an $n\times n$ matrix. 
However, there are minimal $k$-positivity tests which only require testing $n^2$ minors. 
These minimal tests can be related by series of exchanges, and form a family of sub-cluster algebras of the cluster algebra of total positivity tests. 
We give a description of the sub-cluster algebras that give $k$-positivity tests, ways to move between them, and an alternative combinatorial description of many of the tests.
\end{abstract}
\maketitle

\section{Introduction}

A \emph{totally positive} matrix is a matrix in which all minors are positive. Such matrices were originally studied by I. J. Schoenberg in connection with a variation diminishing property~\cite{Schoenberg} and by Gantmacher-Krein due to their nice eigenvalues~\cite{Gantmacher-Krein}.  A \emph{totally nonnegative} matrix is a matrix in which all minors are nonnegative.  Totally positive and totally nonnegative matrices appear in a variety of contexts, including planar networks~\cite{CIM}, canonical bases for quantum groups~\cite{Lusztig1994}, and stochastic processes~\cite{Karlin}.

An $n\times n$ matrix has $\binom{2n}{n}-1$ minors, so it is generally inefficient to test whether a matrix is totally positive by testing all of its minors.  This gives rise to the question of how to test matrices for total positivity as efficiently as possible.  In particular, what are the smallest sets of rational functions in the matrix entries such that positivity of all of these functions ensures total positivity of the matrix?
In~\cite{fomin-zel}, Fomin and Zelevinsky showed that double wiring diagrams give rise to a collection of minimal total positivity tests.  Further, we can obtain a cluster algebra from any double wiring diagram that provides us with additional minimal tests \cite{ClustAlgBook}.

A natural generalization of total positivity is the notion of \emph{$k$-positivity}. Here, we only require that minors of order up to $k$ be positive. We may similarly define \emph{$k$-nonnegativity}, the structure of which is explored in~\cite{knn-bruhat-factorizations}. Tests for $k$-positivity are often discussed in conjunction with the total positivity case \cite{fallat-johnson, pinkus}. These papers give specific classes of tests, but do not explore their combinatorial structure.
In this paper, we generalize both the cluster algebra and double wiring diagram formulations of total positivity tests to $k$-positivity.

We start by giving background and relevant definitions on cluster algebras, total positivity, and double wiring diagrams. Section~\ref{sec:kp} introduces the $k$-positivity cluster algebras and discusses their embedding into the total positivity cluster algebra. We give a construction for how some of these sub-cluster algebras can be augmented with \emph{test variables} to give $k$-positivity tests.  In Section~\ref{sec:k-essential} we define \emph{$k$-essentiality} of minors, which classifies certain minors that are important for testing $k$-positivity, and identify classes of minors which are and are not $k$-essential.
In Section~\ref{sec:DwD}, we explain how $k$-positivity behaves in the context of double wiring diagrams and find a family of sub-cluster algebras that produce $k$-positivity tests. We also give an indexing of this family by Young diagrams, which lets us generate specific tests from them.    

Future work to be pursued includes determining all minors which are $k$-essential (including resolving Conjecture~\ref{conj:all-k-ess}), and further, whether such minors are included in every $k$-positivity test.
Additionally, we would like to determine whether there are sub-cluster algebras that give $k$-positivity tests outside of our known family and find a characterization of all double wiring diagrams that can be modified to give $k$-positivity tests.

\section{Background}
\label{sec:bg}

We start by giving a brief overview of relevant background on cluster algebras.  For more detailed and general discussion, see \cite{ClustAlg1}, \cite{marsh2013lecture}, and \cite{ClustAlgBook}. These definitions are reproduced in a slightly modified form below. 
Throughout this paper, we will be using the notation $[n]$ for the set $\{1,2,\ldots,n\}$ and $[i,j]$ for the set $\{i,i+1,\ldots,j-1,j\}$. 

\begin{definition}
A \emph{quiver} is a directed multigraph with no loops or 2-cycles. The vertices are labeled with elements of $[m]$. A directed edge $(i,j)$ will be denoted $i\to j$. A \emph{quiver mutation} of a quiver $Q$ at vertex $j$ is a process, defined as follows, that produces another quiver $\mu_j(Q)$.
\begin{enumerate}
\item For all pairs of vertices $i,k$ such that $i\to j\to k$, create an arrow $i\to k$.
\item Reverse all arrows adjacent to $j$.
\item Delete a maximal collection of 2-cycles.
\end{enumerate}
\end{definition}
\begin{definition}
Let $\mathcal F = \mathbb C(u_1,\ldots, u_m)$ be the field of rational functions over $\mathbb C$ in $m$ independent variables (this is our \emph{ambient field}). A \emph{labeled seed} of geometric type in $\mathcal F$ is a pair $(\tilde{\B{x}}, Q)$ where $\tilde{\B{x}}=(x_1,\ldots, x_m)$ is an algebraically independent generating set of $\mathcal F$ over $\mathbb{C}$ and $Q$ is a quiver on $m$ vertices such that vertices in $[n]$ are called \emph{mutable} and vertices in $[n+1,m]$ are called \emph{frozen}.
We call $\tilde{\B{x}}$ the labeled \emph{extended cluster} of the seed and $\B{x}=(x_1, \ldots, x_n)$ the \emph{cluster}.  The elements $x_1, \ldots, x_n$ are the \emph{cluster variables} and the remaining elements $x_{n+1}, \ldots, x_m$ are the \emph{frozen variables}.

\end{definition}
\begin{definition}
A \emph{seed mutation} at index $j\in [n]$ satisfies $\mu_j( (\tilde{\B{x}}, Q)) = (\tilde{\B{x}}', \mu_j(Q))$, where $x_i' = x_i$ if $i\neq j$ and $x_j'$ satisfies the \emph{exchange relation}
\[
x_j x_j' = \prod_{i\to j} x_i + \prod_{j\to k} x_k,
\]
where arrows are counted with multiplicity. The right hand side is also referred to as the \emph{exchange polynomial}.  Notice that we allow seed mutations only at mutable vertices, not at frozen ones.
\end{definition}
From here on, we will refer to seed mutations simply as mutations.
\begin{definition}
If two quivers or two seeds are related by a sequence of mutations, we say they are \emph{mutation equivalent}.
For some initial seed $(\B{\tilde{x}}, Q)$, let $\chi$ be the union of all cluster variables over seeds which are mutation equivalent to $(\B{\tilde{x}}, Q)$. Let $R= \mathbb C[x_{n+1}, \ldots, x_m]$. Then the \emph{cluster algebra} of rank $n$ over $R$ associated to this initial seed is $\mathcal A = R[\chi]$.
\end{definition}
\begin{definition}
We consider two clusters \emph{equivalent} if they share the same variables, up to permutation. The \emph{exchange graph} of a cluster algebra is a graph on vertices indexed by equivalence classes of clusters, where there is an edge between two vertices if the clusters corresponding to the vertices are connected by a mutation.
\end{definition}

In addition to our cluster algebra background, we will need some definitions that come from the study of total positivity.
\begin{definition}
For an $m\times n$ matrix $X$ and sets $I\subseteq[m]$, $J\subseteq[n]$, we will let $X_{I,J}$ be the submatrix of $X$ where we take rows indexed by $I$ and columns indexed by $J$. If $|I|=|J|=\ell$, the determinant of this submatrix is $|X_{I,J}|$, and we call this a \emph{minor of order $\ell$} or more simply an \emph{$\ell$-minor}.
\end{definition}
A few types of minors will be particularly important to us.
\begin{definition}
A \emph{solid minor} is a minor with rows indexed by $I=[i,i+\ell-1]$ and columns indexed by $J=[j,j+\ell-1]$.  An \emph{initial minor} is a solid minor where $1\in I\cup J$.  A \emph{corner minor} is a solid minor where the associated submatrix is located at the bottom-left or top-right of the whole matrix.  In other words, either $I=[n-\ell+1,n]$ and $J=[\ell]$, or $I=[\ell]$ and $J=[n-\ell+1,n]$.
\end{definition}
\begin{definition}
A \emph{total positivity test} is a set of expressions in the indeterminants $\{x_{i,j}\}_{1\leq i,j\leq n}$ such that an $n\times n$ matrix $M=(m_{i,j})_{1\leq i,j\leq n}$ is totally positive if and only if evaluating these expressions for $x_{i,j}=m_{i,j}$ yields all positive numbers.
\end{definition}

From \cite{fomin-zel}, we know that the minimal size of a total positivity test is $n^2$, and by Theorem 9 of \cite{fomin-zel}, the following is a total positivity test.
\begin{definition}
The \emph{initial minors test} is the positivity test consisting of all $n^2$ initial minors.
\end{definition}
Double wiring diagrams give us a combinatorial interpretation of total positivity tests of minimal size where all the expressions in the test are minors. We start by recalling the appropriate definitions from \cite{fomin-zel}.
\begin{definition}
A \emph{wiring diagram} consists of a family of $n$ piecewise straight lines, all of the same color, such that each line intersects every other line exactly once. A \emph{double wiring diagram} is two wiring diagrams of different colors which are overlaid. 
We will color our diagrams red and blue, and number the lines such that when reading from bottom to top, the left endpoints of the red lines are in decreasing order, and the left endpoints of the blue lines are in increasing order. We also draw the red wires as thin and the blue wires as thick for ease of reading black and white copies.  Each diagram has $n^2$ \emph{chambers}, see Figure~\ref{fig:lex-min-dwd-3}. A chamber is \emph{bounded} if it is enclosed entirely by wires, and is called \emph{unbounded} otherwise.
\end{definition}
We can label a chamber by the tuple $(r,b)$, where $r$ is the subset of $[n]$ indexing all red strings passing below the chamber, and $b$ is the subset of $[n]$ indexing all blue strings passing below the chamber.
\begin{example}
\label{ex:lex-min-dwd}
Figure~\ref{fig:lex-min-dwd-3} gives double wiring diagram with the chambers labeled appropriately.
\begin{figure}[htp]
\centering
\includegraphics[width=\textwidth]{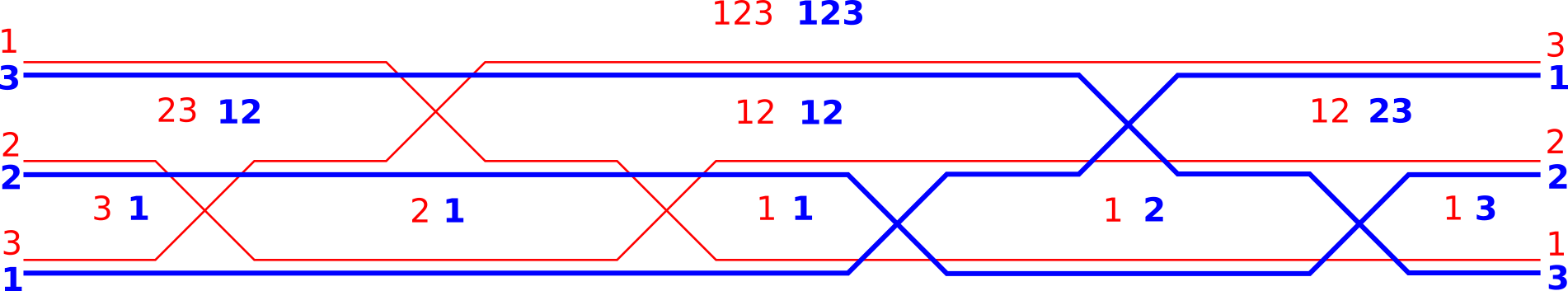}
\caption{A double wiring diagram with wires and chambers labeled. The labels in the chamber correspond to the wires passing underneath.}
\label{fig:lex-min-dwd-3}
\end{figure}
\end{example}
We can associate each chamber with the minor of the corresponding submatrix $\left| X_{r,b} \right|$. With this correspondence, every double wiring diagram gives a total positivity test (Theorem 16 of \cite{fomin-zel}).

\begin{figure}[htp]
\centering
\includegraphics[width=\textwidth]{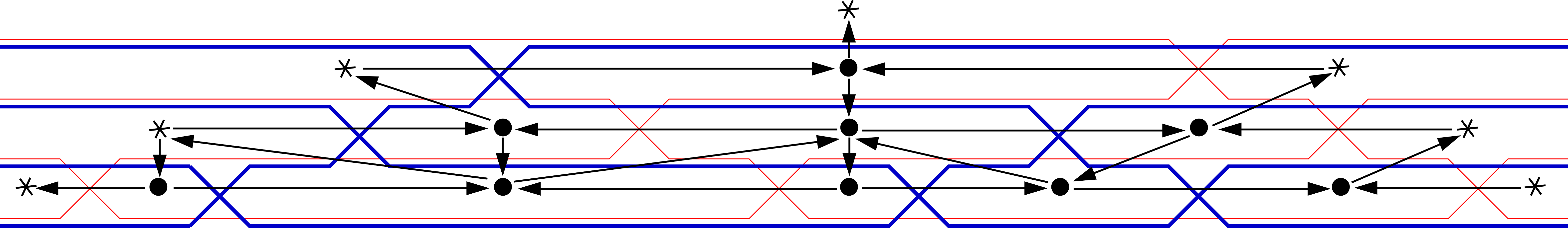}
\caption{An example of $Q(D)$ for a double wiring diagram with $n=4$.}
\label{fig:generic-quiver}
\end{figure}

Additionally, each double wiring diagram can be associated to a quiver, using Definition 2.4.1 of \cite{ClustAlgBook}: 
\begin{definition}
Let $D$ be a double wiring diagram. We construct a quiver $Q(D)$ whose vertices are the chambers. Bounded chambers are mutable vertices and unbounded chambers are frozen vertices. Let $c$ and $c'$ be two chambers, at least one of which is bounded. Then there is an arrow $c\to c'$ in $Q(D)$ if and only if one of the following conditions is met:
\begin{enumerate}
\item the right (resp., left) boundary of $c$ is blue/thick (resp., red/thin), and coincides with the left (resp., right) boundary of $c'$. 
\begin{center}
\includegraphics[width=0.5\textwidth]{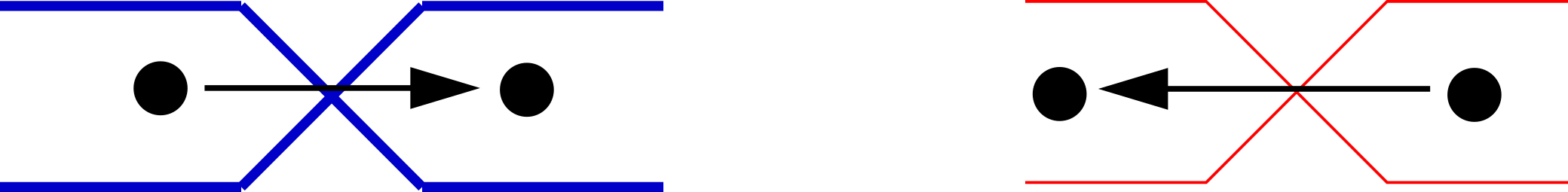}
\end{center}
\item the left boundary of $c'$ is red/thin, the right boundary of $c'$ is blue/thick, and the entire chamber $c'$ lies directly above or directly below $c$.
\begin{center}
\includegraphics[width=0.7\textwidth]{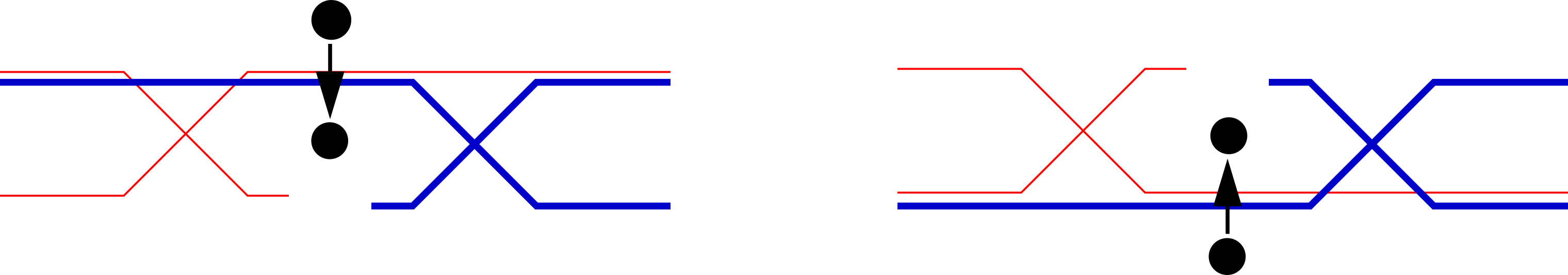}
\end{center}
\item the left boundary of $c$ is blue/thick, the right boundary of $c$ is red/thin, and the entire chamber $c$ lies directly above or directly below $c'$.
\begin{center}
\includegraphics[width=0.7\textwidth]{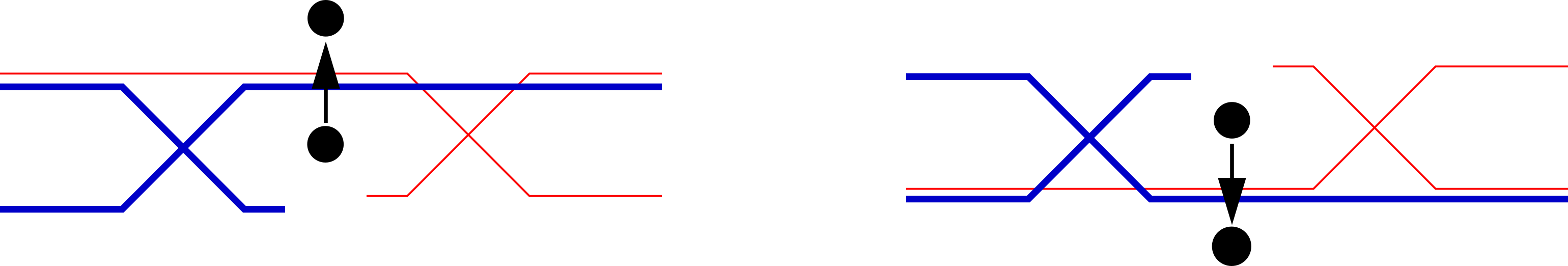}
\end{center}
\item the left (resp., right) boundary of $c'$ is above $c$ and the right (resp., left) boundary of $c$ is below $c'$ and both boundaries are red/thin (resp., blue/thick).
\begin{center}
\includegraphics[width=0.7\textwidth]{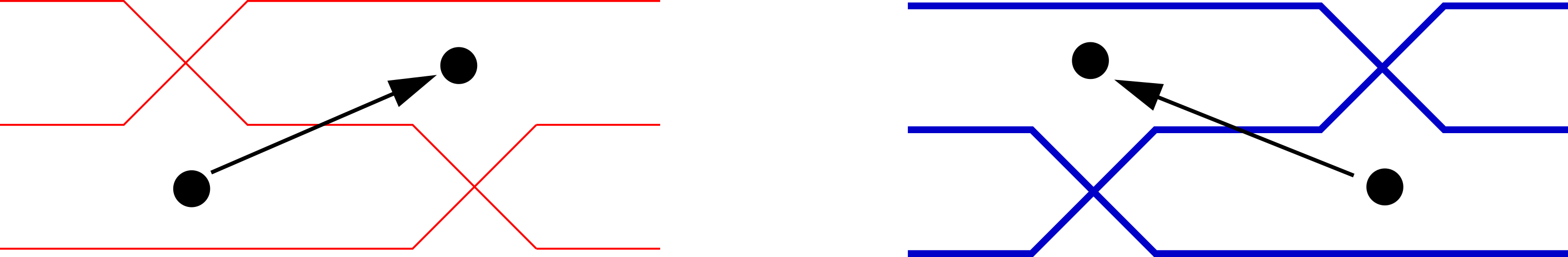}
\end{center}
\item the left (resp., right) boundary of $c$ is above $c'$ and the right (resp., left) boundary of $c'$ is below $c$ and both boundaries are blue/thick (resp., red/thin).
\begin{center}
\includegraphics[width=0.7\textwidth]{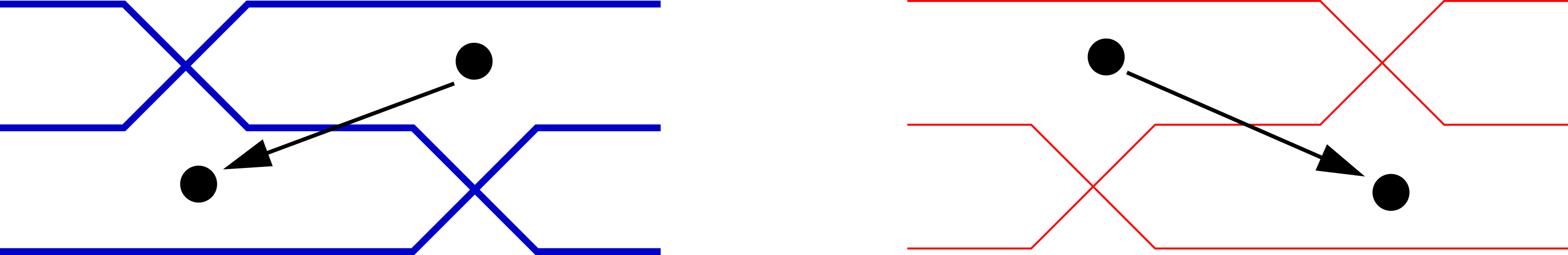}
\end{center}
\end{enumerate}
From each double wiring diagram, we have now shown how to obtain a quiver and how to associate variables to each vertex.  That is, each double wiring diagram gives us a seed.
\end{definition}
Figure~\ref{fig:generic-quiver} shows a generic quiver, illustrating these conditions in context.

\begin{figure}[htp]
\centering
\includegraphics[width=0.85\textwidth]{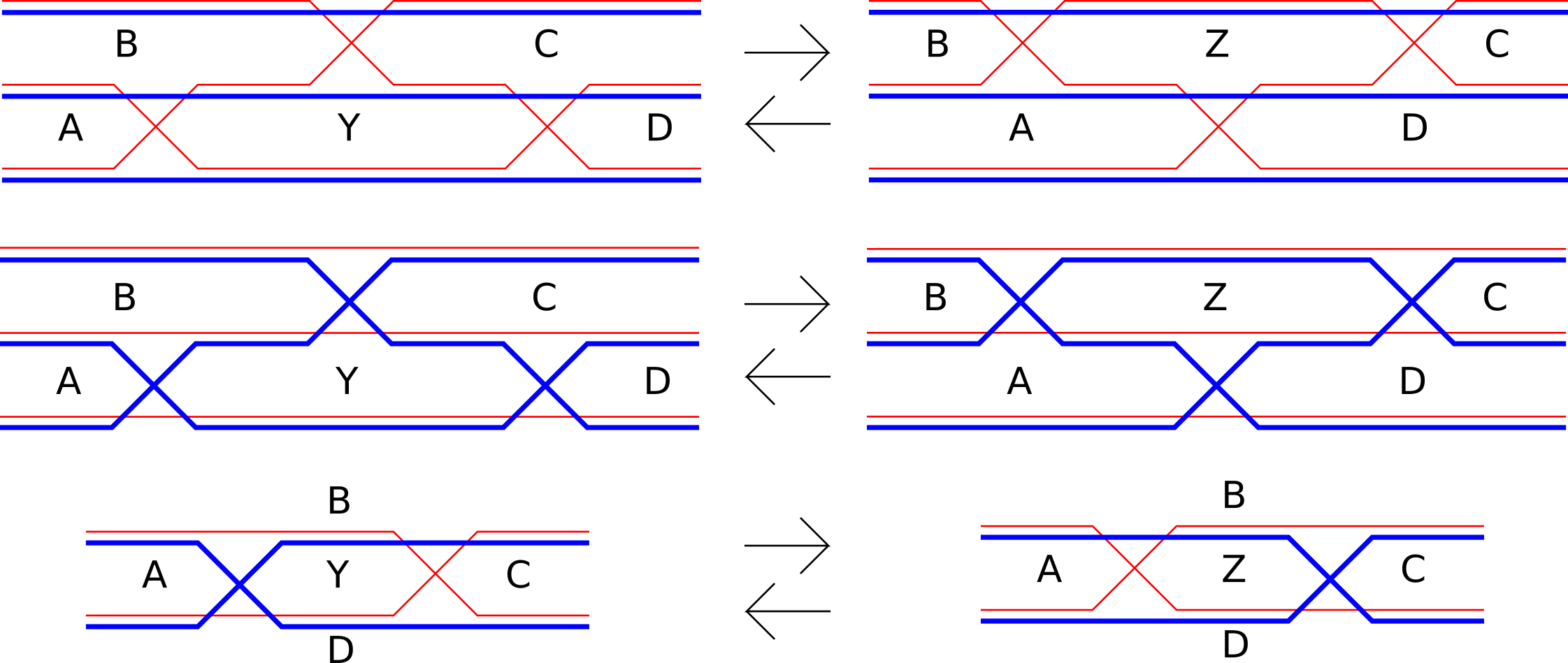}
\caption{The local moves relating double wiring diagrams. The first two are braid moves.}
\label{fig:braid-moves}
\end{figure}

There is also a set of local moves that allow us to transition between any two double wiring diagrams (see \cite{fomin-zel}). These moves are depicted in Figure~\ref{fig:braid-moves}.
The chambers are labeled by the associated minor, and in all cases the exchange relation is $YZ= AC+BD$.  Notice that each of these moves corresponds to mutation at vertex $Y$ in the seed given by the double wiring diagram.

\begin{definition}
The cluster algebra generated by the seed associated to any double wiring diagram is the \emph{total positivity cluster algebra}.
\end{definition}
This cluster algebra has rank $(n-1)^2$ and every seed gives a positivity test.  Notice that in general, this cluster algebra contains seeds that do not correspond to a double wiring diagram.  In fact, there are extended clusters that contain variables that are not matrix minors---other rational functions in the matrix entries can appear. For $n\geq 4$, there are infinitely many such variables, and infinitely many clusters in the total positivity cluster algebra. 

\begin{example}
We will be using the $n=3$ case as a recurring example throughout the next two sections.  For convenience, we'll relabel the entries of our $3\times 3$ as shown below:
\[ M := \begin{bmatrix}
    a & b & c \\
    d & e & f \\
    g & h & j
\end{bmatrix} \]
We can now refer to the $2\times 2$ minors with uppercase letters.  Each uppercase letter will denote the $2\times 2$ minor formed by the rows and columns that do not contain the corresponding lowercase letter.
For example, $A$ is the $2\times 2$ minor obtained from the rows and columns that do not contain $a$, so $A := ej - fh$.
In this case there are only two non-minor extended cluster variables.
These are $K := aA - \det M$ and $L := jJ - \det M$ (see Exercise~1.4.4 of \cite{ClustAlgBook}). 
\end{example}

\section{Generalization to \texorpdfstring{$k$}{k}-Positivity}
\label{sec:kp}

The fact that the total positivity cluster algebra produces total positivity tests relies on two facts:
\begin{enumerate}
\item Every matrix minor appears as a cluster variable.
\item Because each exchange polynomial is subtraction free, when one extended cluster has only positive variables, this means all possible cluster variables must be positive.
\end{enumerate}

If an $n\times n$ matrix is $k$-positive with $k<n$, the variables in an extended cluster of the total positivity cluster algebra are not necessarily all positive.  In fact, if the matrix is not $n$-positive, no extended cluster can have all positive variables.  This poses a problem, as mutation at a vertex corresponding to a positive variable is no longer guaranteed to give us another positive variable.  However, the total positivity cluster algebra leads us to a natural set of sub-cluster algebras that give $k$-positivity tests.

\begin{definition}
Let $(\tilde{\B{x}}, Q)$ be a seed in the total positivity cluster algebra such that every variable in $\B{x}$ is a minor. Designate all vertices in $Q$ corresponding to minors of order greater than $k$ as \emph{dead vertices}, then freeze all vertices in $Q$ which are adjacent to dead vertices.  Allowing mutation only at vertices that are not dead or frozen, we call any quiver that is mutation equivalent to this new quiver a \emph{full $k$-quiver}, and the corresponding seed a \emph{full $k$-seed}.  If we delete the dead vertices as well, we obtain a \emph{$k$-quiver} and \emph{$k$-seed}. 
\end{definition}
This construction resolves the above problem, since now any variable in a seed that is mutation equivalent to such a $k$-seed can be written as a subtraction-free expression in the original $k$-seed's extended cluster variables, which are minors of order at most $k$. These $k$-seeds generate sub-cluster algebras of the total positivity cluster algebra. We have restricted this construction to seeds in the total positivity cluster algebra that have only minors as variables because we don't know whether or not other expressions that show up as cluster variables are required to be positive in $k$-positive matrices.
\begin{example}
Consider the set of quivers in Figure~\ref{fig:n4-k2-set}. Figure~\ref{subfig:n4-quiver} depicts a double wiring diagram and its associated quiver for $n=4$. Here, mutable vertices are represented with $\bullet$ and frozen vertices with $\ast$. Figure~\ref{subfig:n4-full-k2-quiver} shows the full $k$-quiver for $k=2$. Dead vertices are now represented with $\odot$. Figure~\ref{subfig:n4-k2-quiver} shows the $k$-quiver for $k=2$. 
\label{ex:k-quiv}
\begin{figure}[htp]
\begin{subfigure}{\linewidth}
\includegraphics[width=\textwidth]{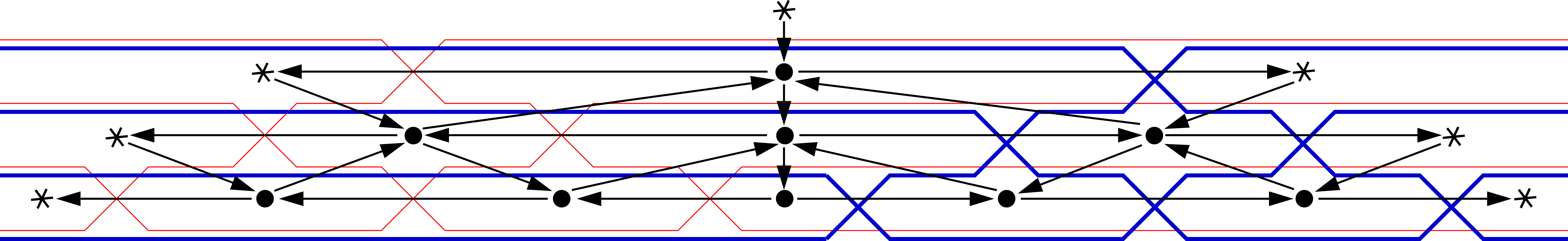}
\caption{The quiver.}
\label{subfig:n4-quiver}
\end{subfigure}

\begin{subfigure}{\linewidth}
\includegraphics[width=\textwidth]{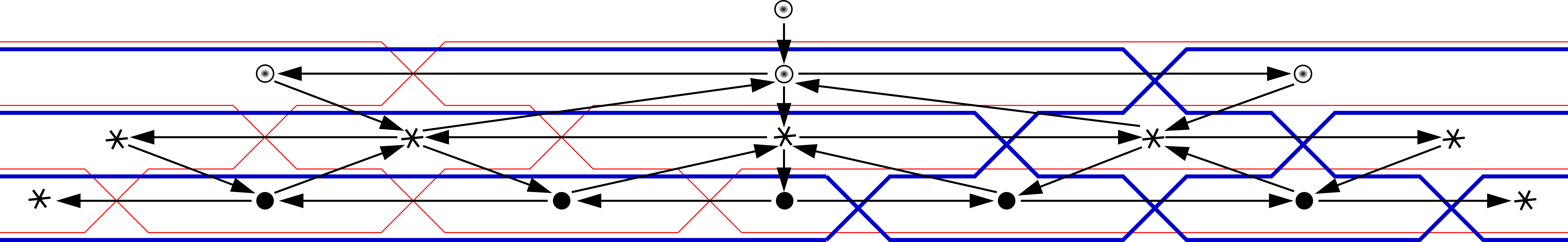}
\caption{The full $k$-quiver.}
\label{subfig:n4-full-k2-quiver}
\end{subfigure}

\begin{subfigure}{\linewidth}
\includegraphics[width=\textwidth]{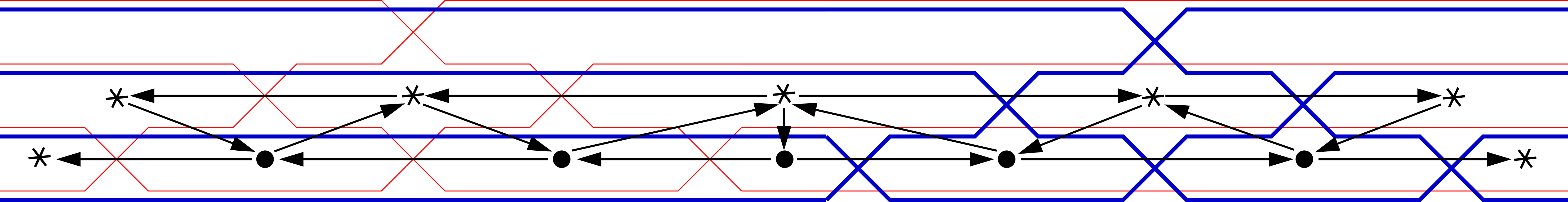}
\caption{The $k$-quiver.}
\label{subfig:n4-k2-quiver}
\end{subfigure}
\caption{The quiver, full $k$-quiver, and $k$-quiver for a particular double wiring diagram with $n=4$, $k=2$. Mutable vertices are depicted by $\bullet$, frozen vertices by $\ast$, and dead vertices by $\odot$.}
\label{fig:n4-k2-set}
\end{figure}
\end{example}
\begin{definition} 
Freezing a vertex in a seed in the total positivity cluster algebra corresponds to deleting all edges corresponding to mutation at that vertex from the exchange graph, and likewise for marking a vertex as dead. The exchange graph for each sub-cluster algebra generated by a $k$-seed is thus a connected component of this new graph, and we will refer to these exchange graphs as \emph{components}.
\end{definition}
\begin{example}
\label{ex:components}
We return to the $n=3$ case, now with $k=2$. For a matrix which is totally positive, $K$ and $L$ (our non-minor variables occurring in tests) must also be positive since they occur in clusters, and hence can be written as subtraction-free rational expressions in the initial minors. For a matrix which is maximally $2$-positive, $K$ and $L$ are also positive as they are both a nonpositive term subtracted from a positive one. 
We can now in this case expand our construction of the sub-cluster algebras to seeds which also contain $K$ or $L$.  In these quivers, we only freeze vertices adjacent to the determinant  and the determinant is the only dead vertex. The exchange graphs for the 8 sub-cluster algebras are depicted in Figure~\ref{fig:TP-comps}. 
The vertices in this figure are labeled by the cluster variables which are mutable in the total positivity cluster algebra, so that the extended cluster contains the listed variables plus $c$, $g$, $C$, and $G$. 

\begin{figure}[thb]
\centering
\begin{tikzpicture}[scale=0.8]
\node(A1) at (7,2) {\tiny$ABfj$};
\node(B1) at (7,4) {\tiny$ABDj$};
\node(C1) at (8,1) {\tiny$Afhj$};
\node(D1) at (8,3) {\tiny$ADhj$};
\node(E1) at (10.5,3.5) {\tiny$ADbh$};
\node(F1) at (12,1.25) {\tiny$Abeh$};
\node(G1) at (10,0) {\tiny$Aefh$};
\node(H1) at (8.5,5.5) {\tiny$ABDL$};
\node(I1) at (10.75,5) {\tiny$ADLb$};
\node(J1) at (12,4) {\tiny$ALbd$};
\node(K1) at (12.5,2.5) {\tiny$Abde$};
\node(L1) at (9.75,4.25) {\tiny$ABLd$};
\node(M1) at (9.25,2.75) {\tiny$ABdf$};
\node(N1) at (10.5,1.75) {\tiny$Adef$};
\node(A2) at (0,2) {\tiny$HJKf$};
\node(B2) at (0,4) {\tiny$FHJK$};
\node(C2) at (1,1) {\tiny$JKfh$};
\node(D2) at (1,3) {\tiny$FJKh$};
\node(E2) at (3.5,3.5) {\tiny$FJbh$};
\node(F2) at (5,1.25) {\tiny$Jbeh$};
\node(G2) at (3,0) {\tiny$Jefh$};
\node(H2) at (1.5,5.5) {\tiny$FHJa$};
\node(I2) at (3.75,5) {\tiny$FJab$};
\node(J2) at (5,4) {\tiny$Jabd$};
\node(K2) at (5.5,2.5) {\tiny$Jbde$};
\node(L2) at (2.75,4.25) {\tiny$HJad$};
\node(M2) at (2.25,2.75) {\tiny$HJdf$};
\node(N2) at (3.5,1.75) {\tiny$Jdef$};
\node(A3) at (3.27,-1.62) {\tiny$BEHj$};
\node(B3) at (5,-4) {\tiny$BEHa$};
\node(C3) at (3.27,-6.38) {\tiny$BHad$};
\node(D3) at (0.48,-5.47) {\tiny$BHdf$};
\node(E3) at (0.48,-2.53) {\tiny$BHfj$};
\node(A4) at (10.27,-1.62) {\tiny$DEFj$};
\node(B4) at (12,-4) {\tiny$DFhj$};
\node(C4) at (10.27,-6.38) {\tiny$DFbh$};
\node(D4) at (7.48,-5.47) {\tiny$DFab$};
\node(E4) at (7.48,-2.53) {\tiny$DEFa$};
\node(A5) at (0.5,-12) {\tiny$BLad$};
\node(B5) at (0.5,-8) {\tiny$BDLa$};
\node(C5) at (4.5,-8) {\tiny$DLab$};
\node(D5) at (4.5,-12) {\tiny$Labd$};
\node(A6) at (7.5,-12) {\tiny$HKfj$};
\node(B6) at (7.5,-8) {\tiny$Kfhj$};
\node(C6) at (11.5,-8) {\tiny$FKhj$};
\node(D6) at (11.5,-12) {\tiny$FHKj$};
\node(A7) at (15,-9) {\tiny$BDEa$};
\node(B7) at (15,-5) {\tiny$BDEj$};
\node(A8) at (15,-2) {\tiny$EFHj$};
\node(B8) at (15,2) {\tiny$EFHa$};
\path[-,font=\large, >=angle 90, line width=0.5mm]
(A1) edge (B1)
(A1) edge (C1)
(B1) edge (D1)
(D1) edge (C1)
(D1) edge (E1)
(E1) edge (F1)
(F1) edge (G1)
(G1) edge (C1)
(B1) edge (H1)
(H1) edge (I1)
(I1) edge (E1)
(I1) edge (J1)
(J1) edge (K1)
(K1) edge (F1)
(K1) edge (N1)
(N1) edge (G1)
(N1) edge (M1)
(M1) edge (L1)
(L1) edge (J1)
(L1) edge (H1)
(M1) edge (A1)
(A2) edge (B2)
(A2) edge (C2)
(B2) edge (D2)
(D2) edge (C2)
(D2) edge (E2)
(E2) edge (F2)
(F2) edge (G2)
(G2) edge (C2)
(B2) edge (H2)
(H2) edge (I2)
(I2) edge (E2)
(I2) edge (J2)
(J2) edge (K2)
(K2) edge (F2)
(K2) edge (N2)
(N2) edge (G2)
(N2) edge (M2)
(M2) edge (L2)
(L2) edge (J2)
(L2) edge (H2)
(M2) edge (A2)
(A3) edge (B3)
(B3) edge (C3)
(C3) edge (D3)
(D3) edge (E3)
(E3) edge (A3)
(A4) edge (B4)
(B4) edge (C4)
(C4) edge (D4)
(D4) edge (E4)
(E4) edge (A4)
(A5) edge (B5)
(B5) edge (C5)
(C5) edge (D5)
(D5) edge (A5)
(A6) edge (B6)
(B6) edge (C6)
(C6) edge (D6)
(D6) edge (A6)
(A7) edge (B7)
(B8) edge (A8);
\end{tikzpicture}
\caption{The components of a 2-positivity test graph derived from the $3\times 3$ exchange graph.}\label{fig:TP-comps}
\end{figure}
\end{example}

As stated in Section~\ref{sec:bg}, we know that the minimal size of a total positivity test is $n^2$.
Also, we can see that the minimal size of a 1-positivity test is $n^2$, as the entries of the matrix are independent variables that must all be positive.
\begin{conjecture}
For any $k$, the minimal size of a $k$-positivity test is $n^2$.
\end{conjecture}
We will be looking for $k$-positivity tests of size $n^2$.

\begin{definition}
A \emph{potential test cluster} is an extended cluster from a $k$-seed with additional rational functions in the matrix entries appended to the cluster to give a set of size $n^2$.
The variables that are in the test cluster and not the extended cluster are called \emph{potential test variables}.  If the potential test cluster gives a $k$-positivity test, it is called a \emph{test cluster} and the additional variables are called \emph{test variables}. These test cluster variables along with the $k$-seed give a \emph{test seed}.
\end{definition}
\begin{example}
\label{ex:tc-k=n}
All test clusters for $k=n$ are extended clusters in the total positivity cluster algebra.
\end{example}
\begin{example}
\label{ex:k-init}
From Theorem~2.3 of~\cite{fallat-johnson-sokal}, we know that the set of all solid $k$-minors and all initial minors of order less than $k$ gives a $k$-positivity test of size $n^2$.  This is the \emph{$k$-initial minors test}.  The $k$-seed from Figure~\ref{subfig:n4-k2-quiver} includes exactly the initial minors of order $\leq k$.  So, this cluster can be augmented to a test cluster by including all the missing solid $k$-minors as test variables.
\end{example}

Not all choices of potential test variables will give a valid test cluster.  Further, not all clusters can be extended to a test cluster, as we shall discuss in Examples~\ref{ex:k-ess} and~\ref{ex:square-comp}. Although we do know which test variables to add to a cluster to obtain a test cluster in specific cases (see Example~\ref{ex:components} and Theorem~\ref{thm:young-dwd-correspondence}), as of now we lack a proof for the general method.
\begin{remark}
Suppose we have a set of potential test variables and we append these variables to every cluster in a cluster algebra generated by a certain $k$-seed.  Proving that a single potential test cluster from this cluster algebra is a $k$-positivity test proves that all the potential test clusters are $k$-positivity tests: we can go between the variables in the extended clusters using subtraction-free rational expressions, and the rest of the variables in the test cluster stay the same. 
\end{remark}

\begin{definition}
Two test seeds from different sub-cluster algebras have a \emph{bridge} between them if they have the same test cluster and there is a quiver mutation connecting them which occurs at a vertex which is frozen in the $k$-quiver. 
\end{definition}
We can think of a bridge as swapping a cluster variable for a test variable.  This allows us to relate different components using test clusters.
\begin{example}
The two largest components in the $n=3,\ k=2$ case (see Figure~\ref{fig:TP-comps}) both generate 2-positivity tests.  The left associahedron contains $(J,a,b,d,c,g,C,G)$, and so appending the test variable $A$ gives the $k$-initial minors test.  The right associahedron contains the extended cluster $(A,f,h,j,c,g,C,G)$, and so appending the test variable $J$ gives the antidiagonal flip of the $k$-initial minors test.  This is also a $k$-positivity test by Theorem~1.4.1 of~\cite{TNNbook}.  There are four bridges between these components, which we get by swapping the roles of $A$ (a cluster variable on the left and test variable on the right) and $J$ (a test variable on the left and cluster variable on the right) (see Figure~\ref{fig:assoc_correspondence}).

\begin{figure}
\begin{tikzpicture}
\node(A) at (12.5,2) {\tiny$ABfj$};
\node(B) at (12.5,4) {\tiny$ABDj$};
\node(C) at (11.5,1) {\tiny$Afhj$};
\node(D) at (11.5,3) {\tiny$ADhj$};
\node(E) at (9,3.5) {\tiny$ADbh$};
\node(F) at (7.5,1.25) {\tiny$Abeh$};
\node(G) at (9.5,0) {\tiny$Aefh$};
\node(H) at (11,5.5) {\tiny$ABDL$};
\node(I) at (8.75,5) {\tiny$ADLb$};
\node(J) at (7.5,4) {\tiny$ALbd$};
\node(K) at (7,2.5) {\tiny$Abde$};
\node(L) at (9.75,4.25) {\tiny$ABLd$};
\node(M) at (10.25,2.75) {\tiny$ABdf$};
\node(N) at (9,1.75) {\tiny$Adef$};
\node(AA) at (0,2) {\tiny$HJKf$};
\node(BB) at (0,4) {\tiny$FHJK$};
\node(CC) at (1,1) {\tiny$JKfh$};
\node(DD) at (1,3) {\tiny$FJKh$};
\node(EE) at (3.5,3.5) {\tiny$FJbh$};
\node(FF) at (5,1.25) {\tiny$Jbeh$};
\node(GG) at (3,0) {\tiny$Jefh$};
\node(HH) at (1.5,5.5) {\tiny$FHJa$};
\node(II) at (3.75,5) {\tiny$FJab$};
\node(JJ) at (5,4) {\tiny$Jabd$};
\node(KK) at (5.5,2.5) {\tiny$Jbde$};
\node(LL) at (2.75,4.25) {\tiny$HJad$};
\node(MM) at (2.25,2.75) {\tiny$HJdf$};
\node(NN) at (3.5,1.75) {\tiny$Jdef$};
\path[-,font=\large, >=angle 90, line width=0.5mm]
(A) edge (B)
(A) edge (C)
(B) edge (D)
(D) edge (C)
(D) edge (E)
(E) edge (F)
(F) edge (G)
(G) edge (C)
(B) edge (H)
(H) edge (I)
(I) edge (E)
(I) edge (J)
(J) edge (K)
(K) edge (F)
(K) edge (N)
(N) edge (G)
(N) edge (M)
(M) edge (L)
(L) edge (J)
(L) edge (H)
(M) edge (A)
(AA) edge (BB)
(AA) edge (CC)
(BB) edge (DD)
(DD) edge (CC)
(DD) edge (EE)
(EE) edge (FF)
(FF) edge (GG)
(GG) edge (CC)
(BB) edge (HH)
(HH) edge (II)
(II) edge (EE)
(II) edge (JJ)
(JJ) edge (KK)
(KK) edge (FF)
(KK) edge (NN)
(NN) edge (GG)
(NN) edge (MM)
(MM) edge (LL)
(LL) edge (JJ)
(LL) edge (HH)
(MM) edge (AA);
\path[dashed,font=\large, >=angle 90, line width=0.2mm]
(F) edge (FF)
(G) edge (GG)
(K) edge (KK)
(N) edge (NN);
\end{tikzpicture}
    \caption{The bridges between the two largest components in the $n=3$, $k=2$ case. The left has test variable $A$ and the right has test variable $J$. There are 4 bridges between these components, which we obtain by matching $(J,d,e,f)$-$(A,d,e,f)$, $(J,e,f,h)$-$(A,e,f,h)$, $(J,b,e,h)$-$(A,b,e,h)$, and $(J,b,d,e)$-$(A,b,d,e)$, 
    i.e. those with the same test cluster (which also includes variables $c$, $g$, $C$, $G$ in all cases).
    }
 \label{fig:assoc_correspondence}
\end{figure}
\end{example}

\begin{remark}
If one sub-cluster algebra provides $k$-positivity tests, then so do any sub-cluster algebras connected by a bridge. This is easy to see because the test cluster that both sub-cluster algebras share is a $k$-positivity test, which tells us that all test clusters in the second sub-cluster algebra are $k$-positivity tests.
\end{remark}

\section{\texorpdfstring{$k$}{k}-essential minors}
\label{sec:k-essential}

To help determine which of the components provide tests, we define the following:
\begin{definition}
A minor $\left| X_{I,J} \right|$ is \emph{k-essential} if $|I|=|J|\leq k$ and there exists a matrix M such that $\left| M_{I,J} \right| \leq 0$, but $ \forall (I',J') \neq (I,J)$, $|I'|=|J'|\leq k$, we have $\left| M_{I',J'} \right| > 0$.
\end{definition}
That is to say, a $k$-essential minor must appear in all $k$-positivity tests consisting only of minors.
\begin{remark}
By the combinatorial proof of Theorem~3.1.10 of \cite{TNNbook} and the discussion following it, all corner minors are $n$-essential. Note that if an $\ell$-minor is $k$-essential, then that minor is $k'$-essential for all $\ell \leq k'\leq k$. Thus, all corner minors of order $\ell$ are $k$-essential for $\ell\leq k$.
\end{remark}

The following two propositions detail additional minors that are $k$-essential in certain cases.

\begin{proposition}
\label{prop:k-essential-2}
Solid 2-minors are 2-essential.
\end{proposition}
\begin{proof}

Let $I=\{i,i+1\}$, $J=\{j,j+1\}$, and consider the matrix
\[ M := \begin{bmatrix} 
    \ddots& \vdots & \vdots & \vdots & \vdots & \vdots & \vdots & \iddots \\
    \cdots& \varepsilon^{-5} & \varepsilon^{-3} & 1 & \varepsilon^3 & \varepsilon^7 & \varepsilon^{10} & \cdots \\
    \cdots& \varepsilon^{-3} & \varepsilon^{-2} & 1 & \varepsilon^2 & \varepsilon^5 & \varepsilon^7 & \cdots \\
    \cdots& 1 & 1 & \varepsilon & 1 & \varepsilon^2 & \varepsilon^3 & \cdots \\
    \cdots& \varepsilon^3 & \varepsilon^2 & 1 & \varepsilon & 1 & 1 & \cdots \\
    \cdots& \varepsilon^7 & \varepsilon^5 & \varepsilon^2 & 1 & \varepsilon^{-2} & \varepsilon^{-3} & \cdots \\
    \cdots& \varepsilon^{10} & \varepsilon^7 & \varepsilon^3 & 1 & \varepsilon^{-3} & \varepsilon^{-5} & \cdots \\
    \iddots & \vdots & \vdots & \vdots & \vdots & \vdots & \vdots & \ddots \\
\end{bmatrix} \]
where $\varepsilon$ is a sufficiently small positive constant.

We will show that the following construction makes all 2-minors except $M_{I, J}$ positive. $M$ is defined so that
\[
M_{I,J} = \begin{bmatrix}
	\varepsilon & 1 \\
    1 & \varepsilon
\end{bmatrix}
\]
and the powers of $\varepsilon$ throughout the rest of the matrix are inductively chosen as follows. 
Fill the rest of rows $i$, $i+1$ and columns $j$, $j+1$ with consecutive increasing powers as shown. We now inductively fill in the rest of the matrix:
\begin{itemize}
\item For $i' > i+1$, $j' > j+1$, let $m_{i',j'} = \varepsilon^{k_{i',j'}}$, where 
\[
k_{i', j'} = \min\{k_{i', t} + k_{s, j'} - k_{s, t} - 1\ |\ i \leq s < i', j \leq t < j'\}.
\]
\item For $i' < i$, $j' > j+1$, let $m_{i',j'} = \varepsilon^{k_{i',j'}}$, where 
\[
k_{i', j'} = \max\{k_{i', t} + k_{s, j'} - k_{s, t} + 1\ |\ i' < s \leq i+1, j \leq t < j'\}.
\]
\item For $i' > i+1$, $j' < j$, let $m_{i',j'} = \varepsilon^{k_{i',j'}}$, where 
\[
k_{i', j'} = \max\{k_{i', t} + k_{s, j'} - k_{s, t} + 1\ |\ i \leq s < i', j' < t \leq j+1\}.
\]
\item For $i' < i$, $j' < j$, let $m_{i',j'} = \varepsilon^{k_{i',j'}}$, where 
\[
k_{i', j'} = \min\{k_{i', t} + k_{s, j'} - k_{s, t} - 1\ |\ i' < s \leq i+1, j < t \leq j+1\}.
\]
\end{itemize}
Now consider $M'$, identical to $M$ but with $M_{I,J}$ replaced by
\[
M'_{I,J} = \begin{bmatrix}
	1 & \varepsilon \\
    \varepsilon & 1
\end{bmatrix}
\]
In this matrix, all solid $2$-minors are positive, because a solid $2$-minor must be entirely in a single quadrant plus the center cross and thus by construction is positive. Applying the $k$-initial minors test, this matrix is then 2-positive.

All 2-minors of $M$ that do not have any entries in $M_{I,J}$ are the same as those in $M'$, and therefore are positive.  All other 2-minors in $M$, except for $|M_{I,J}|$, are positive by construction.  Therefore, $|M_{I,J}|$ is the only non-positive minor of size 2 or less in $M$.
\end{proof}

\begin{proposition}
\label{prop:k-essential-3}
Solid 3-minors are 3-essential.
\end{proposition}
\begin{proof}
This time there are more cases, but the argument is roughly analogous. Let $I=\{i,i+1,i+2\}$, $J=\{j,j+1,j+2\}$ and consider the matrix

\[ M := \begin{bmatrix}
    \ddots & \vdots & \vdots & \vdots & \vdots & \vdots & \iddots \\
    \cdots & \varepsilon^{-2} & 1 & \varepsilon^2 & \varepsilon^8 & \varepsilon^{14} & \cdots \\
    \cdots & 1 & 1+\varepsilon & 1+\varepsilon & \varepsilon^4 & \varepsilon^8 & \cdots\\
    \cdots & \varepsilon^2 & 1+\varepsilon & 1+2\varepsilon & 1+\varepsilon & \varepsilon^2 & \cdots \\
    \cdots & \varepsilon^8 & \varepsilon^4 & 1+\varepsilon &1+\varepsilon & 1 & \cdots\\
    \cdots & \varepsilon^{14} & \varepsilon^8 & \varepsilon^2 & 1 & \varepsilon^{-2} & \cdots \\
    \iddots & \vdots & \vdots & \vdots & \vdots & \vdots & \ddots\\
\end{bmatrix} \]

It is defined so that
\[
M_{I,J} = \begin{bmatrix}
	1+\varepsilon & 1+\varepsilon & \varepsilon^4 \\
    1+\varepsilon & 1+2\varepsilon & 1+\varepsilon \\
    \varepsilon^4 & 1+\varepsilon & 1+\varepsilon
\end{bmatrix}
\]
where $\varepsilon$ is a sufficiently small positive constant, and the powers of $\varepsilon$ throughout the rest of the matrix are inductively chosen as follows. For $j'>j+2$, let 
\begin{align*}
m_{i, j'} &= \varepsilon^{4(j'-j) - 4} &
m_{i+1, j'} &= \varepsilon^{2(j'-j) - 4} &
m_{i+2, j'} &= 1.
\end{align*}
Proceed symmetrically, as seen above in $M$, for the rows in $I$ where $j' < j$, and for columns in $J$ in the regions $i' > i+2$ and $i' < i$. For $i'\in I$ or $j'\in J$, define $k_{i',j'}$ to be the lowest power of $\varepsilon$ found in $m_{i',j'}$. For the four corner regions, as in the case of 2-essentiality, let $m_{i',j'} = \varepsilon^{k_{i',j'}}$, where 
\begin{align*}
k_{i', j'} &= \min\{k_{i', t} + k_{s, j'} - k_{s, t} - 2\ |\ i \leq s < i', j \leq t < j'\}
\end{align*}
for $i' > i+2$, $j' > j+2$, 
\begin{align*}
k_{i', j'} &= \max\{k_{i', t} + k_{s, j'} - k_{s, t} + 2\ |\ i \leq s < i', j' < t \leq j+2\}
\end{align*}
for $i' > i+2$, $j' < j$, and symmetrically for the other two regions.

Now, consider $M'$, where $M_{I,J}$ has been replaced with

\[M'_{I,J} := \begin{bmatrix}
    1 & \varepsilon & \varepsilon^4 \\
    \varepsilon & 1 & \varepsilon \\
    \varepsilon^4 & \varepsilon & 1 \\
\end{bmatrix} \]

We apply

\begin{lemma}\label{lemma:e-matrix}
Any 2-positive matrix $M$ whose entries are all powers of $\varepsilon$ is totally positive when $\varepsilon$ is made sufficiently small.
\end{lemma}

\begin{proof}
Since any submatrix is also a 2-positive matrix whose entries are all power of $\varepsilon$, it suffices to show that the determinant of the whole matrix is positive. This determinant is the sum of entries of the form $\pm \Pi m_{i, \sigma(i)}$. It suffices to show that compared to the diagonal term, all other terms are smaller by a factor of at least $\varepsilon$. This fact follows from the positivity of all 2-minors: since $\sigma$ can be decomposed into transpositions, any term can be created from the diagonal term by repeatedly replacing $m_{i,j} m_{i',j'}$ with $m_{i,j'} m_{i',j}$ for some $i<i'$, $j<j'$. Positivity of the minor $|M_{\{i,i'\},\{j,j'\}}|$ gives that $\varepsilon m_{i,j} m_{i',j'} \geq m_{i',j} m_{i,j'}$.
\end{proof}

As in the case of 2-essentiality, all solid $2$-minors in $M'$ are positive, because a solid $2$-minor must lie entirely in a single quadrant plus the center cross and thus by construction is positive. Applying the $k$-initial minors test, this matrix is then 2-positive. Then, by Lemma~\ref{lemma:e-matrix}, $M'$ is totally positive, and any minors it shares with $M$ (namely those with no entries in $M_{I,J}$) are therefore positive. By construction, all 2-minors in $M$ with an entry in $M_{I,J}$ are also positive.

In fact, the method used to prove Lemma~\ref{lemma:e-matrix} works for all 3-minors except those which contain a 2-minor not satisfying $\varepsilon m_{i,j} m_{i',j'} \geq m_{i,j'} m_{i',j}$.  $M$ has exactly two such minors, namely $|M_{\{i,i+1\},\{j,j+1\}}|$ and $|M_{\{i+1,i+2\},\{j+1,j+2\}}|$.

Now consider $M_{\{i_1, i_2, i_3\},\{j_1, j_2, j_3\}}$ containing exactly one of $M_{\{i,i+1\},\{j,j+1\}}$ and $M_{\{i+1,i+2\},\{j+1,j+2\}}$. Consider any submatrix $M_{\{s,s'\},\{t,t'\}}$ of $M$. For $M_{\{i,i+1\},\{j,j+1\}}$ and $M_{\{i+1,i+2\},\{j+1,j+2\}}$, we can check the following inequality holds:
\begin{align}
\label{eq:half-ineq}
m_{s,t} m_{s',t'} - m_{s',t} m_{s,t'} \geq \frac{1}{2}\varepsilon m_{s,t} m_{s',t'}.
\end{align}
This inequality also holds for all other choices of $M_{\{s,s'\},\{t,t'\}}$, because it follows from the inequality in the previous paragraph.

Suppose $i_1 = i$, $i_2 = i+1$, $j_1 = j$, $j_2 = j+1$. We then also have the following inequalities:
\begin{align}
\label{eq:e2-ineq}
\varepsilon^2 m_{i_1, j_1}m_{i_2, j_2}m_{i_3, j_3}\geq m_{i_1, j_1}m_{i_2, j_3}m_{i_3, j_2}\\
\label{eq:e2-ineq2}
\varepsilon^2 m_{i_1, j_1}m_{i_2, j_2}m_{i_3, j_3}\geq m_{i_1, j_3}m_{i_2, j_2}m_{i_3, j_1}
\end{align}
These conditions are guaranteed by the inductive construction, since $\varepsilon^2 m_{i,j} m_{i',j'} \geq m_{i,j'} m_{i',j}$ unless $\{i,i'\} \subset I$, $\{j,j'\} \subset J$.  This gives us the following:
\begin{align*}
|M_{\{i_1, i_2, i_3\},\{j_1, j_2, j_3\}}| \geq&\ m_{i_1, j_1}m_{i_2, j_2}m_{i_3, j_3} - m_{i_1, j_2}m_{i_2, j_1}m_{i_3, j_3} - m_{i_1, j_1}m_{i_2, j_3}m_{i_3, j_2} -\\
&\ m_{i_1, j_3}m_{i_2, j_2}m_{i_3, j_1}\\
\geq&\ \frac{1}{2}\varepsilon m_{i_1, j_1}m_{i_2, j_2}m_{i_3, j_3} - m_{i_1, j_1}m_{i_2, j_3}m_{i_3, j_2} -
m_{i_1, j_3}m_{i_2, j_2}m_{i_3, j_1}\\
&\ \text{by Inequality~\ref{eq:half-ineq}}\\
\geq&\ \left(\frac{1}{2}\varepsilon - 2\varepsilon^2\right) m_{i_1, j_1}m_{i_2, j_2}m_{i_3, j_3}\text{ by Inequalities~\ref{eq:e2-ineq} and~\ref{eq:e2-ineq2}} \\
>&\ 0. 
\end{align*}

We have analogous inequalities for any other 3-minor containing exactly one of the 2-minors in question.  Thus the only 3-minor of $M$ not necessarily positive is one containing both such 2-minors, $|M_{I,J}|$ itself. It can be checked that this minor is in fact negative, thus showing that it is 3-essential.
\end{proof}

\begin{example}
\label{ex:k-ess}
Returning to the $n=3$, $k=2$ case, the previous proposition tells us that $A$ and $J$ are $2$-essential.
From Figure~\ref{fig:TP-comps}, we can see that the pentagonal components and single-edge components in the $n=3$, $k=2$ case are all missing both $A$ and $J$ from their extended clusters. The extended clusters are of size 8, and in this case $n^2=9$.  All of the extended cluster variables in these components are minors, and $A$ and $J$ are 2-essential, so these components cannot give 2-positivity tests of size $n^2$ which are composed entirely of minors.
\end{example}

Based on the above results and the prevalence of solid $k$-minors in our other $k$-positivity tests, we propose the following:
\begin{conjecture}
\label{conj:all-k-ess}
Solid $k$-minors are $k$-essential.
\end{conjecture}

The $k=1$ case is trivial, as there is exactly one 1-positivity test containing only matrix minors: the test consisting of all $n^2$ elements of the matrix. Explicitly, we can let $x_{i,j} = -1, x_{i', j'} = 1$ for $(i,j) \neq (i',j')$. The cases of $k=2$, $k=3$ are proven in Propositions~\ref{prop:k-essential-2} and~\ref{prop:k-essential-3}, respectively.

Unfortunately, the technique used in Propositions~\ref{prop:k-essential-2} and~\ref{prop:k-essential-3} fails in the general case.  This is because the central minors called $M_{I,J}$ in the above propositions were generated from a maximally 1-nonnegative $2\times2$ matrix and a maximally 2-nonnegative $3\times3$ matrix, each consisting only of 1's and 0's. Theorem~2.2 of \cite{brualdi-kirkland} shows that no such maximally $(k-1)$-nonnegative $k\times k$ matrices exist for $k>3$. This has made extrapolating to the general case rather difficult.

While our definition of $k$-essentiality tells us about which minors must be present in $k$-positivity tests composed entirely of minors, it would also be beneficial to know more about which minors must be present in every test composed of cluster variables.
\begin{example}
\label{ex:square-comp}
Consider the matrix
\[\begin{bmatrix}
    \varepsilon & 1 & \varepsilon^2 \\
    1 & \varepsilon & 1 \\
    \varepsilon^2 & 1 & \varepsilon^{-2}
\end{bmatrix} \]
for some small positive constant $\varepsilon$. All the minors of orders 1 and 2 are positive, except for $J$. In addition, the non-minors $K$ and $L$ are also positive. Thus, the positivity of $J$ is not implied by the positivity of any other cluster variables in the $n=3$, $k=2$ case. This means $J$ must appear in every 2-positivity test which only uses cluster variables from the $n=3$, $k=2$ cluster algebra, regardless of whether the test contains non-minors.
Using the antidiagonal flip of this matrix, we can see that the same holds for $A$.
From Figure~\ref{fig:TP-comps}, we can see that all of the extended clusters in the pentagonal, square, and single-edge components are missing both $A$ and $J$.  As in Example~\ref{ex:k-ess}, the extended clusters are all of size 8 and $n^2=9$.  So these components cannot give 2-positivity tests of size $n^2$ using cluster variables.
\end{example}

It is also useful to know which minors aren't $k$-essential.
\begin{proposition}
If $|I|=|J|<k$ and $|X_{I,J}|$ is not a corner minor, then $|X_{I,J}|$ is not $k$-essential.
\end{proposition}
\begin{proof}
Suppose for a matrix $M$ that $|M_{I,J}|\leq 0$. We show that there exists some other minor of order at most $k$ which is also nonpositive. Pick indices $i\in I$, $j\in J$, $i'\notin I$, $j'\notin J$ such that either $i<i'$ and $j<j'$ or $i>i'$ and $j>j'$. This is possible since $|I|=|J|<n$ and the minor isn't a corner minor.

For ease of reading, we will omit brackets around sets containing one or two elements in the following.

Recall Lewis Carroll's identity: if $M$ is an $n\times n$ square matrix and $M_A^B$ is $M$ with the rows indexed by $A$ and columns indexed by $B$ removed, then 
$$\det(M_a^b)\det(M_{a'}^{b'})-\det(M_a^{b'})\det(M_{a'}^b) = \det(M)\det(M_{a,a'}^{b,b'})$$
if $1\leq a<a'\leq n$ and $1\leq b<b'\leq n$.
Using this identity on the matrix $M_{I\cup i', J\cup j'}$ gives
\[
|M_{I,J}|\cdot |M_{(I\cup i')\setminus i, (J\cup j')\setminus j}| = |M_{I,(J\cup j')\setminus j}|\cdot |M_{(I\cup i')\setminus i, J}| + |M_{I\cup i', J\cup j'}|\cdot |M_{I\setminus i, J\setminus j}|.
\]
By our initial assumption, $|M_{I,J}|\leq 0$. If $|M_{(I\cup i')\setminus i, (J\cup j')\setminus j}|\leq 0$ then we're done; otherwise the left-hand side is $\leq 0$. Thus at least one summand on the right-hand side must also be $\leq 0$, which means at least one of $|M_{I,(J\cup j')\setminus j}|$, $|M_{(I\cup i')\setminus i, J}|$, $|M_{I\cup i', J\cup j'}|$, $|M_{I\setminus i, J\setminus j}|$, all of which are minors of order at most $k$, is $\leq 0$. Since this holds for all $M$, $|X_{I,J}|$ is not $k$-essential.
\end{proof}


\section{Double Wiring Diagrams}
\label{sec:DwD}

We now return to double wiring diagrams. These will give us a more combinatorial way to think about $k$-positivity tests, and can be used to find different components giving $k$-positivity tests.

To describe a double wiring diagram, it is sufficient to describe the relative positions of all of the crossings. We can think of a diagram as having $n$ tracks numbered from bottom to top, where the chambers in track $i$ have $|r|=|b|=i$ and each crossing occurs in one of the first $n-1$ tracks.  We label a red crossing in the $i^{\text{th}}$ track as $e_i$, and a blue crossing in the $i^{\text{th}}$ track as $f_i$. 
With this notation, a sequence of crossings describing a double wiring diagram is a reduced word for the element $(w_0,w_0)$ of the Coxeter group $S_n\times S_n$, where $w_0$ is the order-reversing permutation (the longest word), see \cite{fomin-zel}. We now define some useful groupings of crossings. Let $r_i = e_{n-i}\cdots e_2e_1$ for $1\leq i \leq n-1$, and let $b_i = f_{1} f_2\cdots f_{n-i}$ for $1\leq i \leq n-1$.
For convenience, when $i\notin [n-1]$ we define $r_i$ and $b_i$ to be empty, containing no crossings. Generally, $r_i$ looks like a diagonal chain of red crossings going down and to the right, starting in the $(n-i)^\text{th}$ track and ending in the first track. Similarly, $b_i$ looks like a a diagonal chain of blue crossings going up and to the right, starting in the first track and ending in the $(n-i)^\text{th}$ track.

\begin{example}
Suppose $n=4$. Then the set of red/thin crossings on the left is $r_1$ and the set of blue/thick crossings on the right is $b_2$.  
\begin{center}
\includegraphics[width=0.5\textwidth]{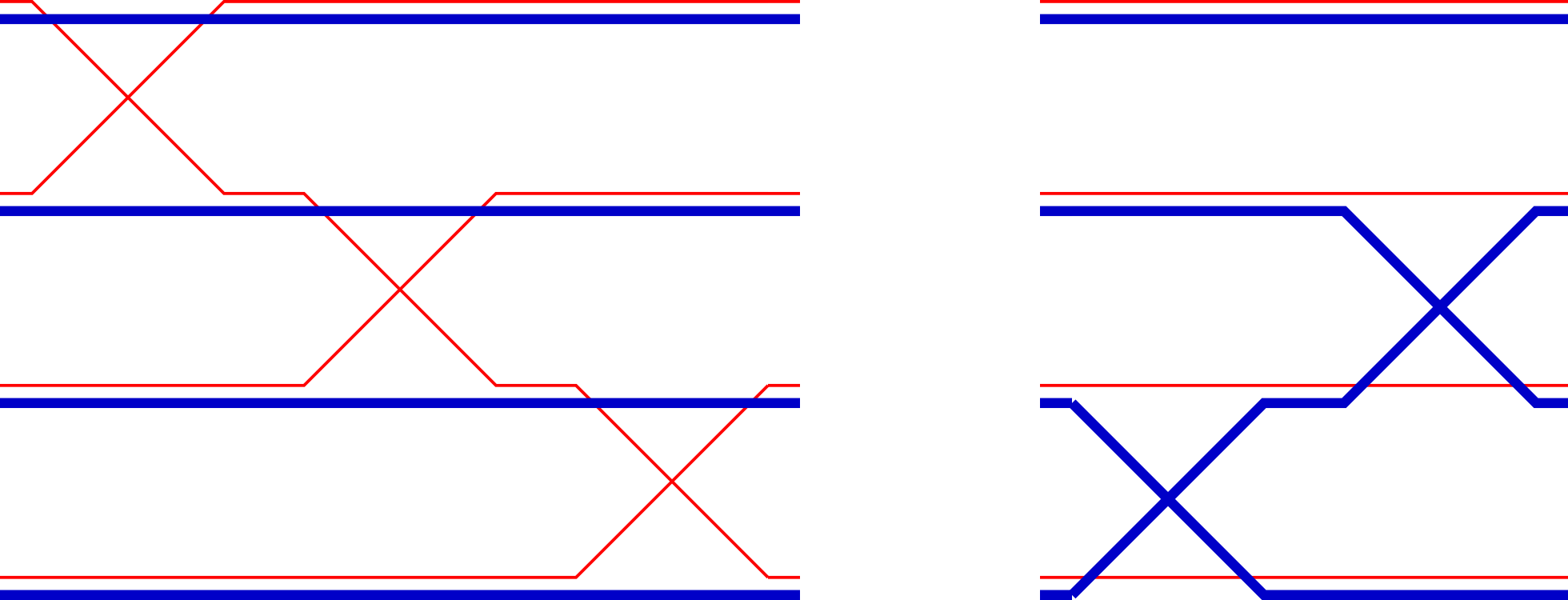}
\end{center}
\end{example}

\begin{definition}
The \emph{lexicographically minimal diagram} is the word $r_{n-1}\cdots r_{1} b_{1} \cdots b_{n-1}$. 
The \emph{lexicographically maximal diagram} is the word $b_{1} \cdots b_{n-1} r_{n-1} \cdots r_{1}$.
\end{definition}

\begin{example}
The lexicographically minimal diagram for $n=3$ appears in Figure~\ref{fig:lex-min-dwd-3} and for $n=4$ appears in Figure~\ref{subfig:n4-quiver}.
\end{example}

\begin{figure}[htp]
\includegraphics[width=\textwidth]{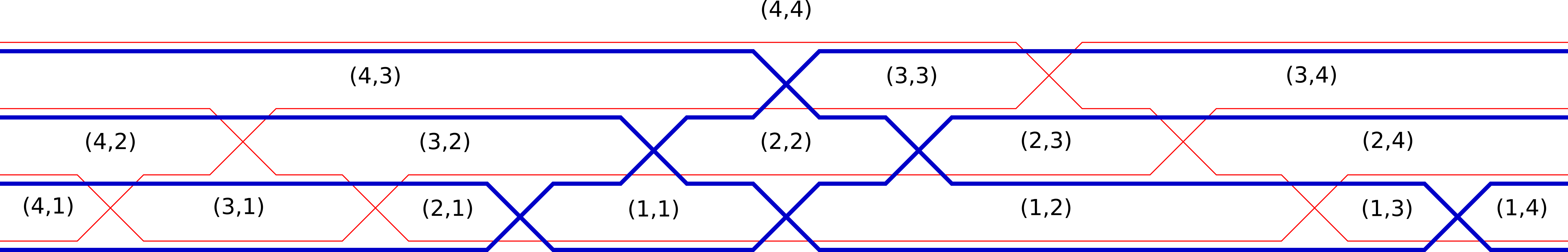}
\caption{The chamber labelling for diagram $r_3r_2b_1b_2r_1b_3$ with $n=4$.}
\label{fig:chamber-index}
\end{figure}

In a diagram given by an interleaving of $b_1\cdots b_{n-1}$ and $r_{n-1}\cdots r_1$, track $t$ always has $2(n-t) +1$ chambers.  We will label the chambers from left to right as $(n,t), (n-1,t),\ldots, (t,t), \ldots, (t,n-1), (t,n)$, as in Figure~\ref{fig:chamber-index}. We 
use this labelling in the following proposition.

\begin{proposition}
\label{prop:lex-min-dwd}
The minor associated to chamber $(i,j)$ in the lexicographically minimal diagram is 
\[
\{|X_{[i-m,i],[j-m,j]}|\}_{1\leq i,j\leq n}
\]
where $m=\min\{i-1,j-1\}$.  The arrows in the quiver are $(i,j) \to (i+1,j)$, $(i,j)\to (i,j+1)$, and $(i+1,j+1)\to (i,j)$ for all $i,j\in[n-1]$.  The vertices $(n,j)$ and $(i,n)$ are frozen for all $i,j\in[n]$.

The minor associated to chamber $(i,j)$ in the lexicographically maximal diagram has variables
\[
\{|X_{[n-j+1, n-j+m+1], [n-i+1, n-i+m+1]}|\}_{1\leq i,j\leq n}
\]
where $m=\min\{i-1, j-1\}$. The arrows in the quiver are $(i+1,j) \to (i,j)$, $(i,j+1)\to (i,j)$, and $(i,j)\to (i+1,j+1)$ for all $i,j\in [1,n-1]$.  The vertices $(n,j)$ and $(i,n)$ are frozen for all $i,j\in[n]$.
\end{proposition}
\begin{proof}
First consider the lexicographically minimal case. 
At the left of our diagram we have the minors $|X_{[n-t+1,n],[t]}|$ in the $t^\text{th}$ track for $t<n$. We inductively proceed by moving left to right, crossing the red groupings.  
Crossing from the left to right of $r_s$, only minors in the first through $n-s^\text{th}$ tracks change. 
The crossings in $r_s$ bring the red wire labelled $s$ to the bottom, and so between $r_s$ and $r_{s-1}$, the red wires are, from bottom to top, $s, s+1, \ldots, n, s-1, s-2, \ldots, 1$ (see Figure~\ref{fig:ri-wire-labels}). The chamber minor between $r_s$ and $r_{s-1}$ in the $t^{\text{th}}$ track is $|X_{[s,s+t-1], [t]}|$ for $1<s<n$ and $t\leq n-s$. 

\begin{figure}[htp]
\centering
\includegraphics[width=0.8\textwidth]{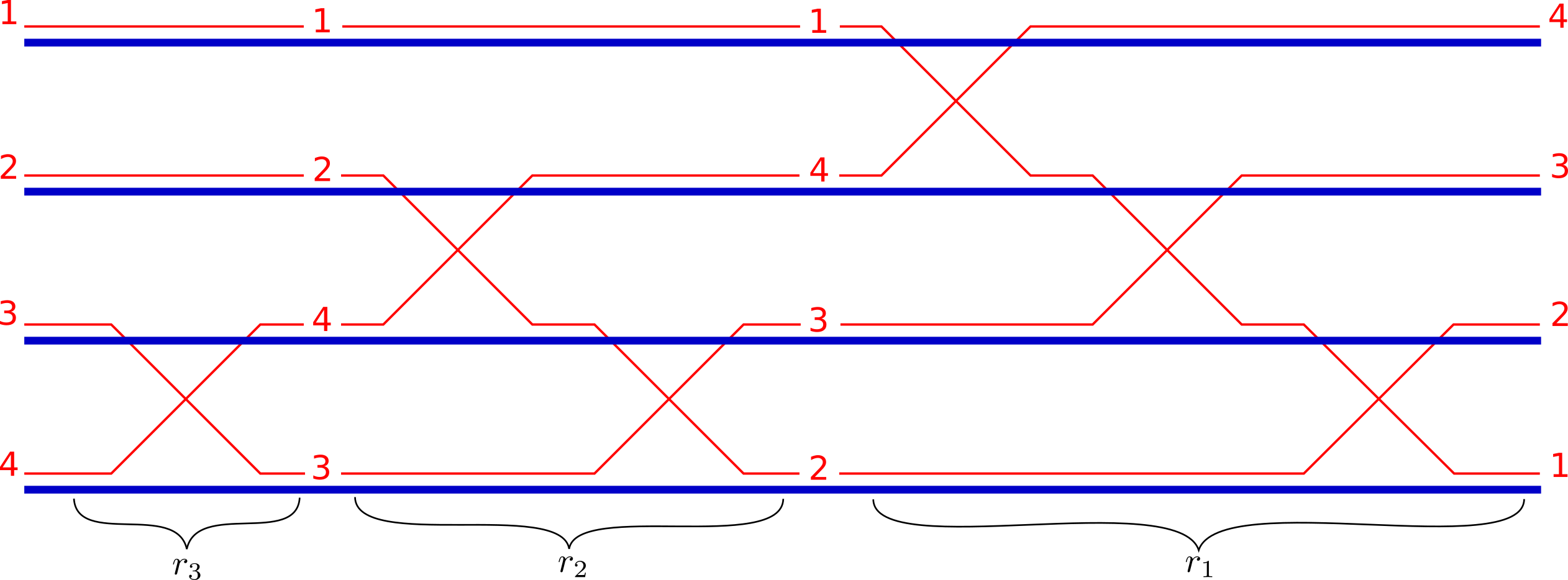}
\caption{Part of the lexicographically minimal wiring diagram showing $r_3, r_2$, and $r_1$. Note that $r_3$ brings the red/thin wire labeled $3$ to the bottom of the diagram and the red/thin wires between $r_3$ and $r_2$ are labelled $3, 4, 2, 1$. We can see that $r_2$ and $r_1$ act similarly.}
\label{fig:ri-wire-labels}
\end{figure}

The chamber minors between $r_1$ and $b_1$ are $|X_{[t], [t]}|$ for $t\leq n-1$.  Now crossing the blue groupings, we see that between $b_{s-1}$ and $b_s$, the blue wires are, from bottom to top, $s, s+1, \ldots, n, s-1, s-2, \ldots, 1$. The chamber minor between $b_{s-1}$ and $b_s$ in the $t^{\text{th}}$ track is $|X_{[t],[s,s+t-1]}|$ for $1<s<n$ and $t\leq n-s$.  On the right of the diagram we have the minors $|X_{[t],[n-t+1,n]}|$ in the $t^\text{th}$ track for $t<n$.  Finally, at the top of our diagram, we have the determinant of the whole matrix.  These are exactly the variables $\{|X_{[i-m,i],[j-m,j]}|\}_{1\leq i,j\leq n}$ where $m=\min\{i-1,j-1\}$. In fact, labeling the vertex corresponding to the variable $|X_{[i-m,i],[j-m,j]}|$ as $(i,j)$, we can read across the $t^\text{th}$ track the vertices $(n,t), (n-1,t),\ldots,(t+1,t),(t,t),(t,t+1),\ldots,(t,n)$.

The description of the arrows can be verified using the definition of $Q(D)$. The arrows $(i,j)\to (i+1, j)$ come from condition 1 (red) and condition 4 (blue). 
The arrows $(i,j) \to (i,j+1)$ come from condition 1 (blue) and condition 4 (red). 
The arrows $(i+1, j+1) \to (i,j)$ come from conditions 2 and 5. Figure~\ref{fig:arrow-index} depicts the $n=4$ case.
\begin{figure}[htp]
\includegraphics[width=\textwidth]{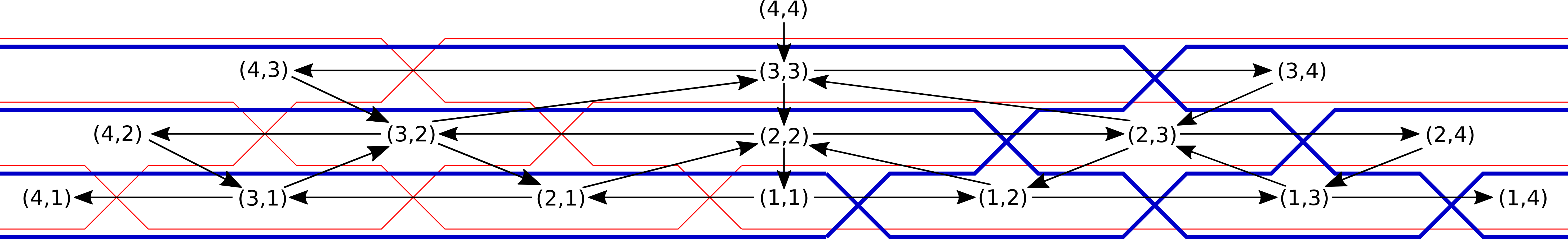}
\caption{The lexicographically minimal diagram for $n=4$ with chambers labelled and quiver arrows shown.}
\label{fig:arrow-index}
\end{figure}

The proof for the lexicographically maximal diagram proceeds in an analogous way.
\end{proof}
Note that the lexicographically minimal diagram gives the initial minors test as defined in Section~\ref{sec:bg}.  The lexicographically maximal diagram gives an antidiagonal flip of this test.  Because of this, we will name the seeds produced by these two diagrams.
\begin{definition}
The seed given by the lexicographically minimal diagram is the \emph{initial minors seed}.
The seed given by the lexicographically maximal diagram is the \emph{antidiagonal initial minors seed}.
\end{definition}

\begin{definition}
The $k$-seed we obtain from the lexicographically minimal diagram has exactly the initial minors of size at most $k$ as its cluster variables.  Thus, we call this $k$-seed the \emph{$k$-initial minors seed}.  Similarly, the $k$-seed obtained from the lexicographically maximal diagram is the \emph{antidiagonal $k$-initial minors seed}.
\end{definition}
\begin{remark}
The sub-quiver induced by mutable vertices of the $2$-initial minors quiver is an orientation of the Dynkin diagram $A_{2n-3}$.  This is also the case for the antidiagonal $2$-initial minors quiver. From Theorem~5.1.3 of \cite{marsh2013lecture} and the discussion in Chapter 6, the component arising from these $2$-quivers is the corresponding associahedron of Cartan type $A_{2n-3}$.  In the $n=3,\ k=2$ example, these components are the two largest components (see Figure~\ref{fig:TP-comps}).
\end{remark}

By definition, a chamber in track $i$ corresponds to a minor of order $i$. Therefore the dead vertices correspond to chambers above track $k$, and no chamber below track $k$ is ever dead or frozen. We note that we can now refer to the braid moves from Figure~\ref{fig:braid-moves} as $e_\ell e_{\ell+1} e_\ell \leftrightarrow e_{\ell+1}e_\ell e_{\ell+1}$ and $f_\ell f_{\ell+1} f_\ell \leftrightarrow f_{\ell+1} f_\ell f_{\ell+1}$.  These braid moves have exchange relations which use minors of orders $\ell$ and $\ell+1$. The other local move from Figure~\ref{fig:braid-moves}, $e_\ell f_\ell \leftrightarrow f_\ell e_\ell$, uses minors of orders $\ell -1$, $\ell$, and $\ell+1$. This means the disallowed local moves are of the form $e_k f_k \leftrightarrow f_k e_k$, $e_k e_{k+1} e_k \leftrightarrow e_{k+1} e_k e_{k+1}$, and $f_k f_{k+1} f_k \leftrightarrow f_{k+1} f_k f_{k+1}$. Thus if a chamber in track $k$ can have a local move applied, the corresponding vertex is frozen. Because of this, we will more generally refer to chambers in track $k$ as frozen.

\begin{definition}
A \emph{path} between a pair of test seeds is a sequence of mutations that takes us from first seed in the pair to second seed such that every mutation yields a seed which can be augmented with test variables to form a test cluster.
\end{definition}
To help us describe such paths, we give a construction for getting a double wiring diagram from a Young diagram. Specifically, let $Y$ be a Young diagram which fits in an $(n-1)\times (n-1)$ square. Now construct the double wiring diagram $D(Y)$ as follows:
\begin{enumerate}
\item Start with the word $b_1 b_2\cdots b_{n-1}$.
\item Let $\ell_k$ be the number of boxes in the $k^{\text{th}}$ row of $Y$.
\item For $k\in [n-1]$, insert $r_k$ between $b_{\ell_k}$ and $b_{\ell_k+1}$.  If there are multiple $r_s$'s between some $b_t$ and $b_{t+1}$, arrange the $r_s$'s in decreasing order from left to right.
\end{enumerate}
The result is an interleaving of the words $b_1\cdots b_{n-1}$ and $r_{n-1}\cdots r_1$. 

\begin{example}
From the Young diagram $Y$ depicted below, we get the word $r_3 b_1 r_2 b_2 b_3 r_1$. This is the double wiring diagram depicted in Figure~\ref{fig:generic-quiver}.
\[
Y = \yng(3,1)
\]
\end{example}

We would like to describe the tests given by these Young diagrams.

Recall that in such interleavings of $b_1\cdots b_{n-1}$ and $r_{n-1}\cdots r_1$, we can label the chambers in track $t$ from left to right as $(n,t),(n-1,t),\ldots, (t,t),\ldots, (t,n-1),(t,n)$. In the following proposition, we associate a minor to these vertices.

\begin{proposition}
\label{prop:yd-minors}
Let $Y$ be a Young diagram and $y_d$ be the number of boxes on the $d^{\text{th}}$ diagonal, where the diagonals are labeled as follows:
\begin{center}
\begin{ytableau}
0 & 1 & 2 & 3 & \hdots\\
-1 & 0 & 1 & 2\\
-2 & -1 & 0 & 1\\
\vdots
\end{ytableau}
\end{center}
Let $\alpha_{i,j}:=\min\{y_{j-i},n-i,n-j\}$.  Then the minor associated to vertex $(i,j)$ of $D(Y)$ is
\[
m_{ij} = 
\begin{cases}
|X_{[1+\alpha_{i,j}, i+\alpha_{i,j}], [1+j-i+\alpha_{i,j}, j+\alpha_{i,j}]}|    & i\leq j, \\
|X_{[1+i-j + \alpha_{i,j}, i+\alpha_{i,j}], [1+\alpha_{i,j}, j+\alpha_{i,j}]}|    & i>j.
\end{cases}
\]
\end{proposition}

\begin{proof}
We proceed by induction on the number of boxes in the diagram. The base case is Proposition~\ref{prop:lex-min-dwd}, where $Y=\emptyset$ is the lex minimal diagram, and $Y$ the full $(n-1)\times(n-1)$ square is the lex maximal diagram. Assume the statement holds for diagrams with $m$ boxes.
Recall also from the proof of Proposition~\ref{prop:lex-min-dwd} that between $r_i$ and $r_{i-1}$ the red wires are, from bottom to top, $i,i+1,\ldots,n,i-1,i-2,\ldots,1$ and that between $b_{i-1}$ and $b_i$ the blue wires are, from bottom to top, $i,i+1,\ldots,n,i-1,i-2,\ldots,1$. 
This still applies even in the interleaved case. 
Now add an $\ell^{\text{th}}$ box to the $k^{\text{th}}$ row, where $k$ is a row such that this is a valid addition. 
The new box is added to the $(\ell-k)^\text{th}$ diagonal, which now has $\min(k,\ell)$ boxes.
In particular, row $k-1$ must have had at least $\ell$ boxes, and so this addition changes the word from $\cdots r_k b_{\ell} \cdots $ to $\cdots b_{\ell} r_k \cdots$. The chambers which change are in tracks $\min(n-k,n-\ell)$ and lower, since these are the tracks in which red and blue crossings are being swapped. 
Originally, the chamber in track $t\leq \min(n-k,n-\ell)$ between $r_k$ and $b_\ell$ was $([k,k+t-1] ,  [\ell,\ell+t-1])$.  After the swap it becomes $([k+1,k+t],  [\ell+1,\ell+t])$. 
This is illustrated in Figure~\ref{fig:ri-bi-swap-chamber}.
\begin{figure}[htp]
\begin{subfigure}{0.47\linewidth}
\includegraphics[width=\textwidth]{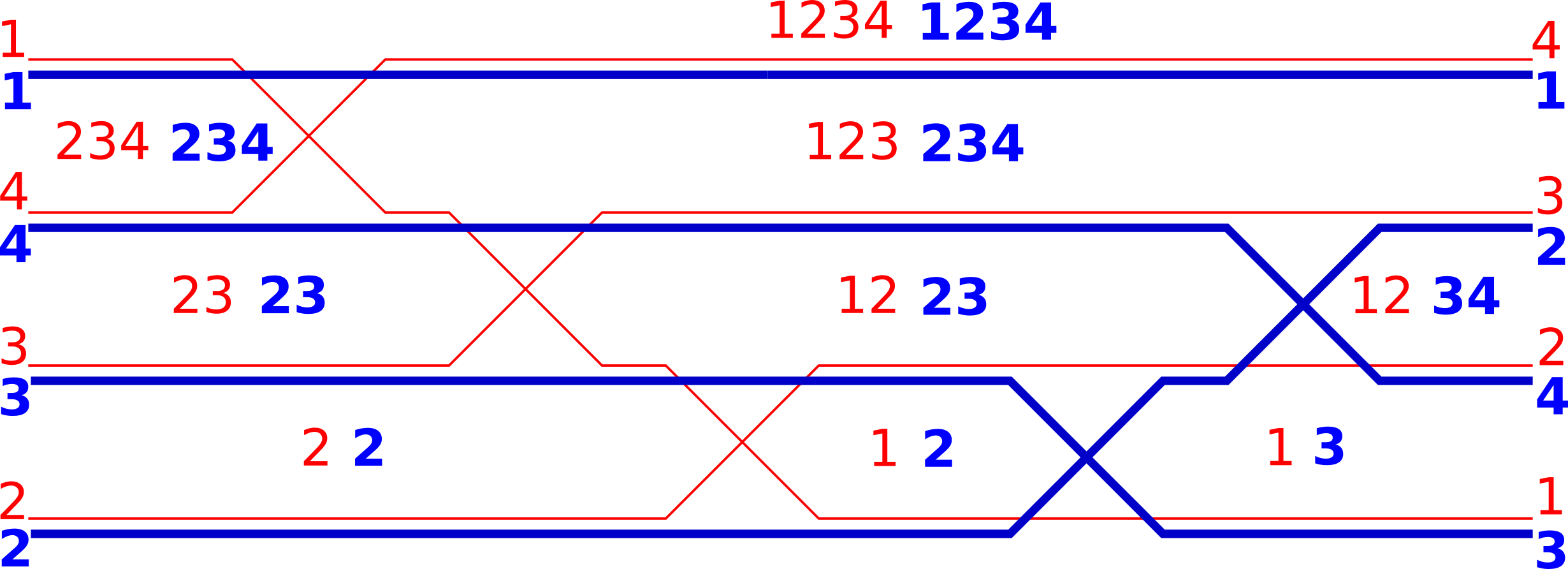}
\caption{$\cdots r_1b_2\cdots $}
\end{subfigure}\quad
\begin{subfigure}{0.47\linewidth}
\includegraphics[width=\textwidth]{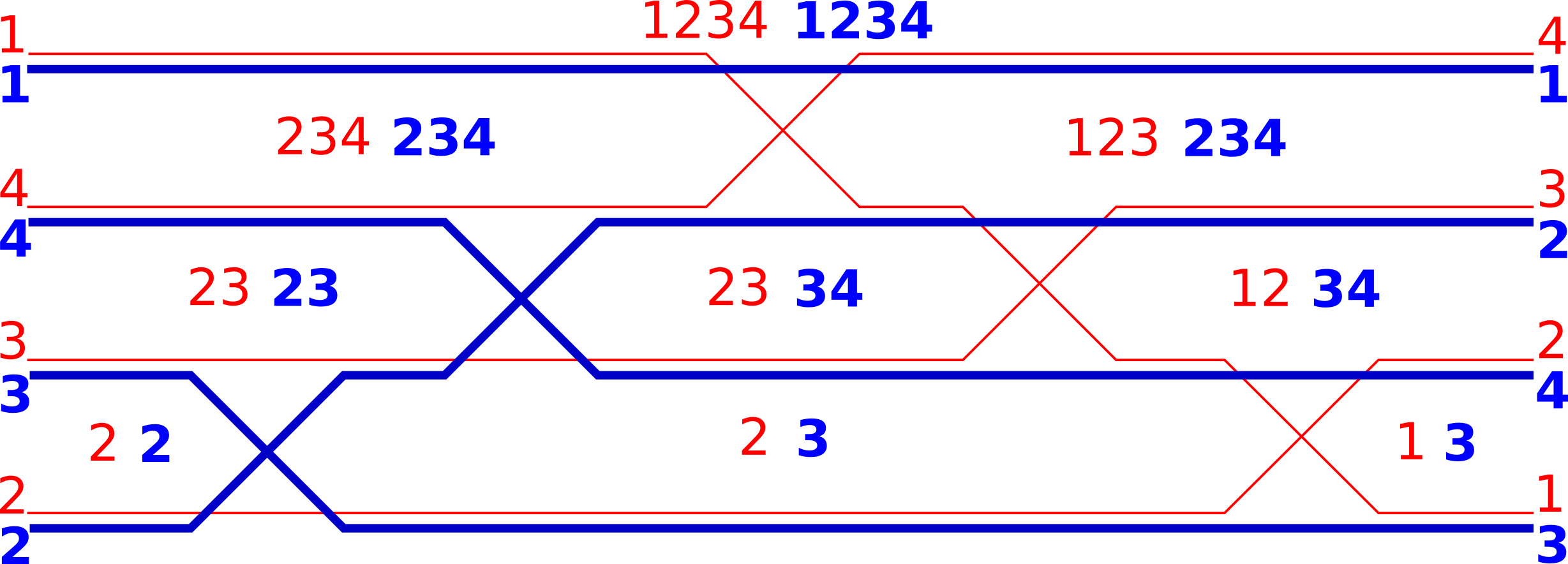}
\caption{$\cdots b_2r_1\cdots$}
\end{subfigure}
\caption{
Swapping $r_1$ and $b_2$ in the $n=4$ case.  Note that the minors corresponding to the middle chambers of tracks 1 and 2 change as specified above.
}
\label{fig:ri-bi-swap-chamber}
\end{figure}

Now we must determine the vertices of $D(Y)$ corresponding to the chambers above that changed. For $t\leq \min(n-k,n-\ell)$, there are $n-t-k$ red crossings to the left of $r_k$ on track $t$, and $k-1$ to the right. Similarly, there are $n-t-\ell$ blue crossings to the right of $b_\ell$ on track $t$, and $\ell-1$ to the left.
So in total, since we also must count $r_k$ and $b_\ell$ themselves, there are $n-t-k+\ell$ crossings on track $t$ to the left of the chamber under consideration, and $n-t-\ell+k$ to the right.
This means there are $n-t-k+\ell$ chambers on track $t$ to the left of the one under consideration, and $n-t-\ell+k$ to the right. The chamber in track $t$ that changed is thus in position
\begin{align*}
(i,j):&=\left( \max(n-\#\textrm{ chambers to left}, t), \ \max(n-\#\textrm{ chambers to right}, t)\right) \\
&=\big( \max(t+k-\ell, t), \ \max(t+\ell-k, t)\big).
\end{align*}
Since $j-i = \ell-k$, we get that the new $y_{j-i}=y_{\ell-k}=\min(k,\ell)$, and so the new $\alpha_{i,j}$ is
\[
\min(k, \ell, n-(t+\max(k-\ell,0)), n-(t+\max(\ell-k,0))) =\min(k,\ell)
\]
under our assumption that $t\leq \min(n-k,n-\ell)$. We can check that these new $\alpha_{i,j}$'s give the correct formulas.  For all other vertices $(i,j)$ the chamber is unchanged and $j-i\neq\ell-k$, so $\alpha_{i,j}$ is also unchanged.
\end{proof}

This Young diagram construction gives us a convenient way to describe certain paths, as the following theorem shows.
\begin{theorem}
\label{thm:young-dwd-correspondence}
Suppose we have a sequence of Young diagrams $Y_0, \ldots, Y_{(n-1)^2}$ such that $Y_0$ is the empty diagram, $Y_{(n-1)^2}$ is the $(n-1)\times (n-1)$ square, and $Y_i$ differs from $Y_{i-1}$ by the addition of a single box. 
This sequence gives a valid path between the lex minimal and lex maximal test seeds by for each added box, mutating the chambers in between the corresponding swapped groupings, working from track 1 upwards.
The $k$-positivity test is formed by disregarding all chambers above the $k^{\text{th}}$ track and adding in the remaining solid $k$-minors, giving a test of size $n^2$.
\end{theorem}
\begin{proof}
The proof that this is an allowed sequence of local moves between the starting and ending double wiring diagrams comes from the first part of the proof of Proposition~\ref{prop:yd-minors}. If we have a Young diagram and add an $\ell\nth$ box to the $i\nth$ row, there must have been at least $\ell$ boxes in row $i-1$ and so the double wiring diagram goes from $\cdots r_i b_\ell \cdots$ to $\cdots b_\ell r_i \cdots$, i.e. from $\cdots e_{n-i}e_{n-i-1}\cdots e_1 f_1 \cdots f_{n-\ell}\cdots$ to $\cdots f_1\cdots f_{n-\ell}e_{n-i}e_{n-i-1}\cdots e_1$. This swap can be formed by applying the local move to swap $e_1$ and $f_1$, then $e_2$ and $f_2$, and so on until we have swapped $e_{\min(n-i, n-\ell)}$ and $f_{\min(n-i,n-\ell)}$. The remaining crossings don't interact with each other and can be slid freely. 
See Figure~\ref{fig:swap-ri-bi} for an illustration.

\begin{figure}[htp]
\begin{subfigure}{0.43\linewidth}
\includegraphics[width=\textwidth]{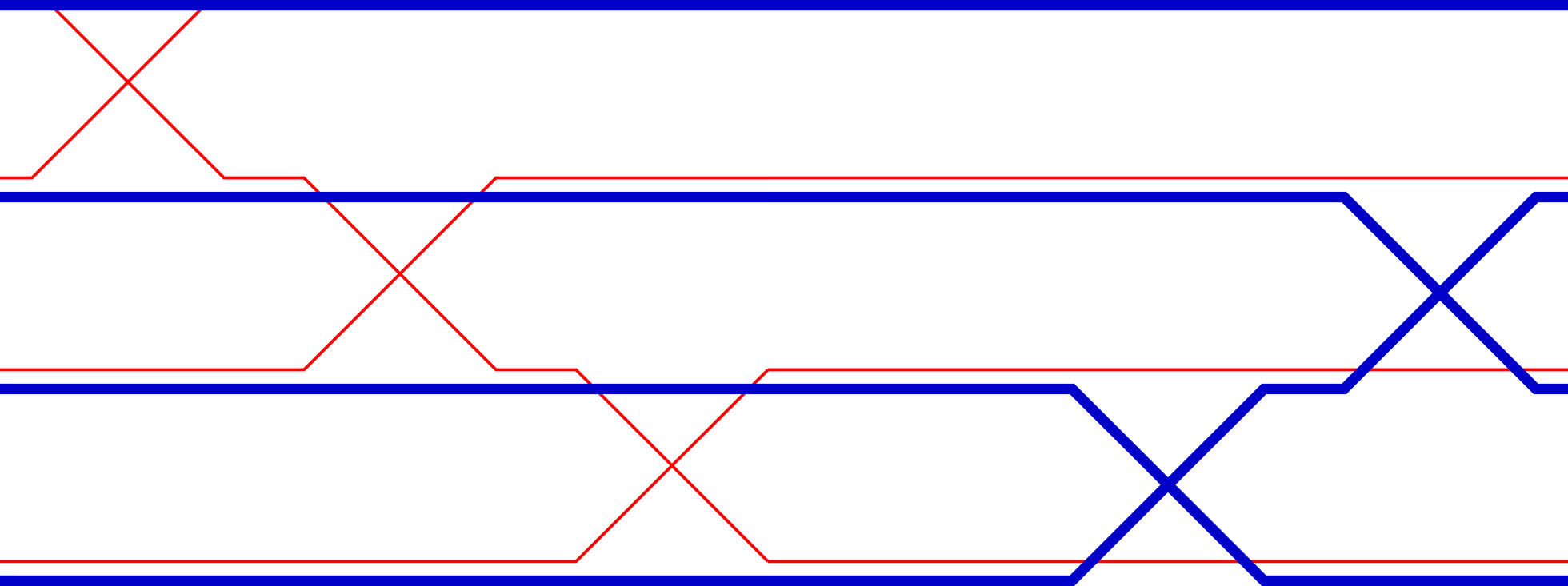}
\caption{Before swapping $r_1$ and $b_2$.}
\label{subfig:swap-init}
\end{subfigure}\qquad
\begin{subfigure}{0.43\linewidth}
\includegraphics[width=\textwidth]{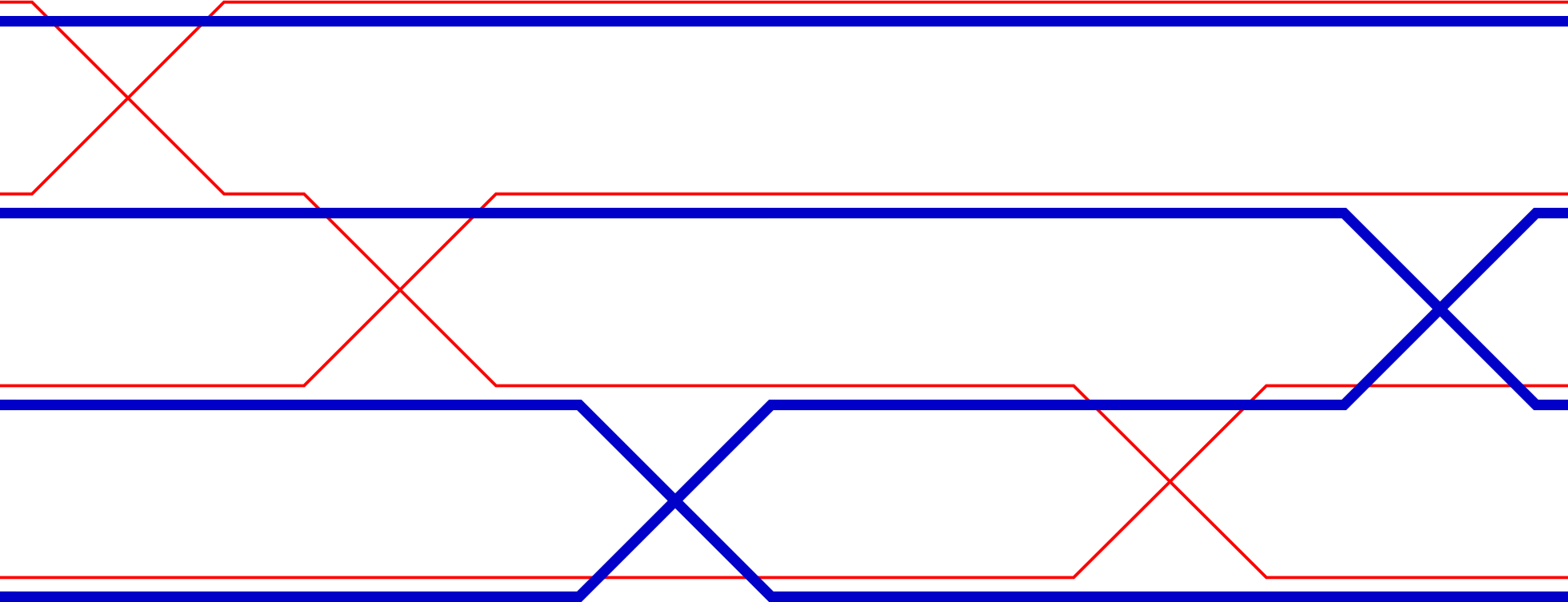}
\caption{Swapping $e_1$ and $f_1$.}
\label{subfig:swap1}
\end{subfigure}
\par\bigskip

\begin{subfigure}{0.43\linewidth}
\includegraphics[width=\textwidth]{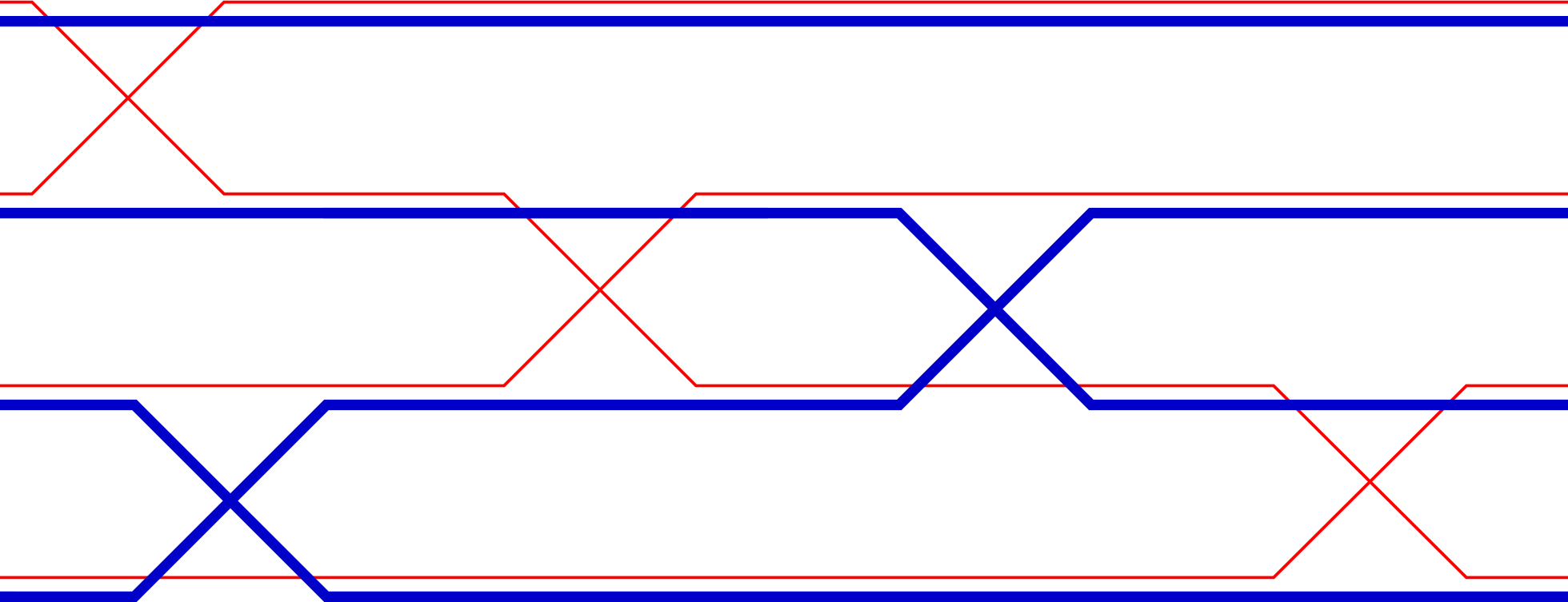}
\caption{Free sliding to prepare for the next swap.}
\label{subfig:swap-between}
\end{subfigure}\qquad
\begin{subfigure}{0.43\linewidth}
\includegraphics[width=\textwidth]{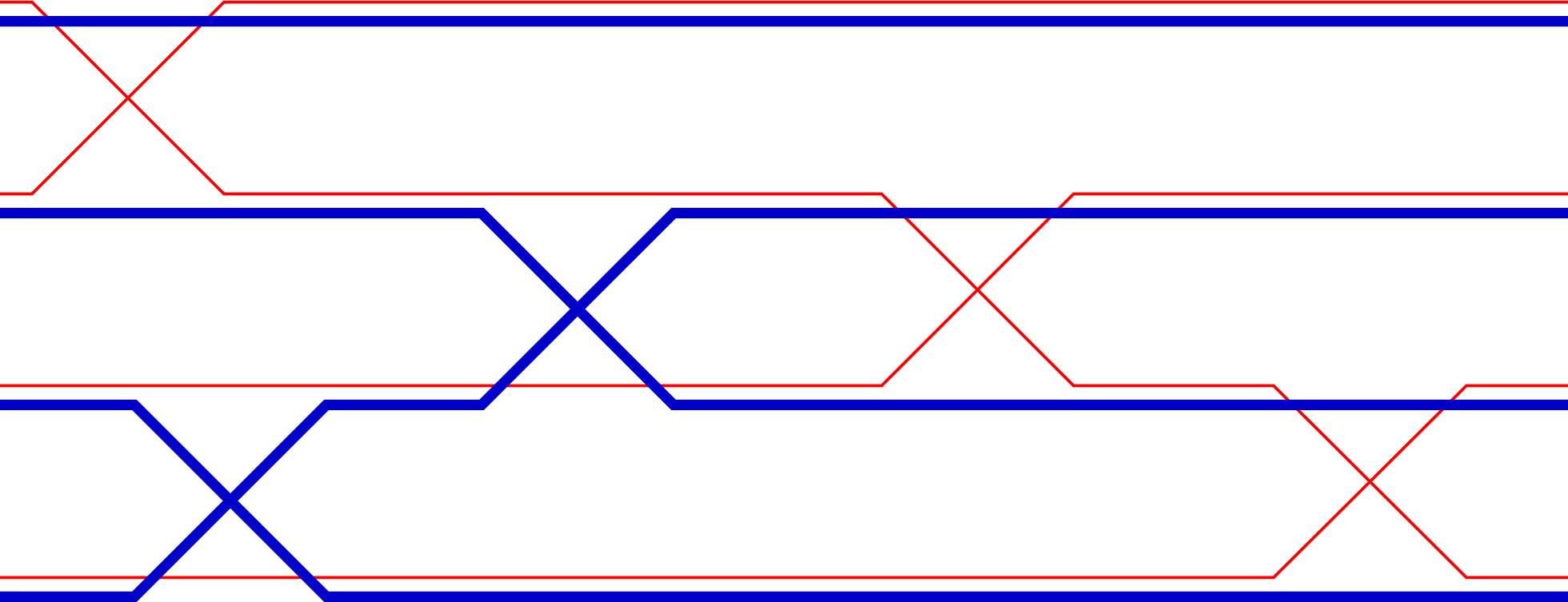}
\caption{Swapping $e_2$ and $f_2$. }
\label{subfig:swap2}
\end{subfigure}%
\par\bigskip

\begin{subfigure}{0.43\linewidth}
\includegraphics[width=\textwidth]{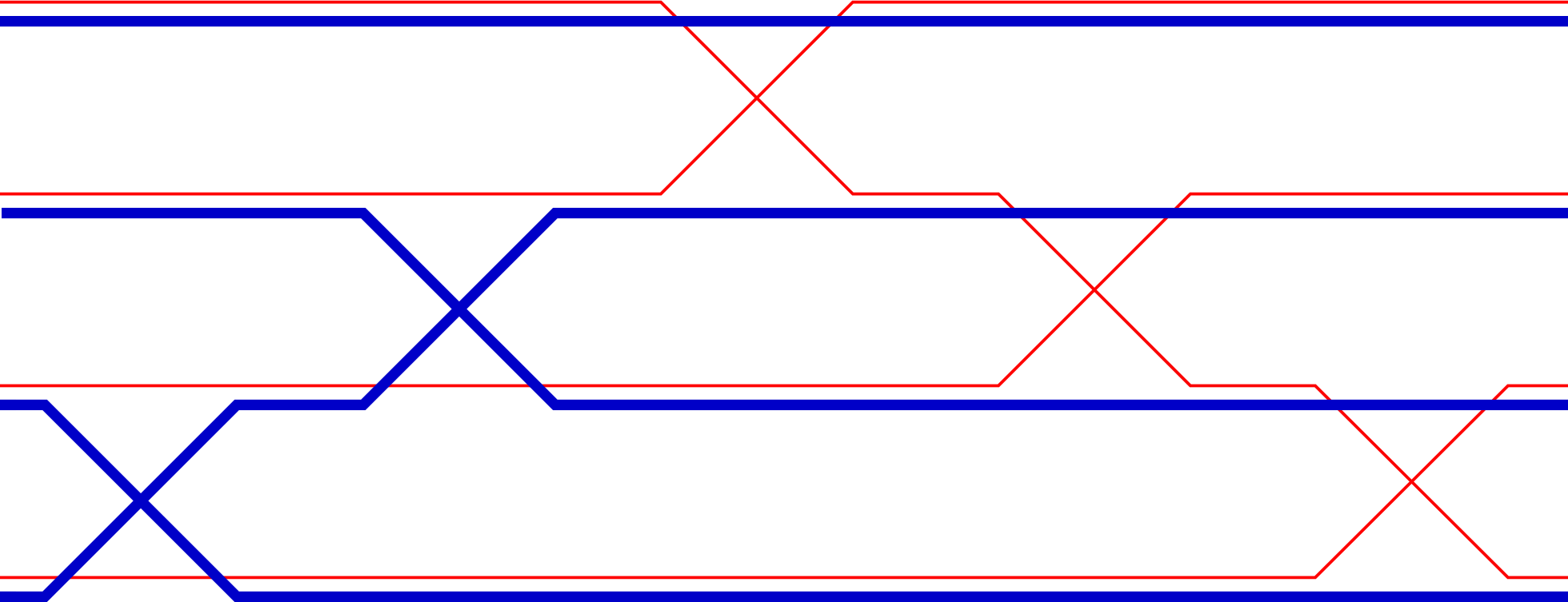}
\caption{Sliding $e_3$ freely to the right.}
\label{subfig:swap-end}
\end{subfigure}
\caption{The stages of swapping $r_1$ and $b_2$ in the $n=4$ case.}
\label{fig:swap-ri-bi}
\end{figure}

Now we confirm that $k$-positivity holds. $D(Y_0)$ gives the $k$-initial minors test. Working inductively, suppose we have diagram $Y_m$ and add an $\ell\nth$ box to the $i\nth$ row.  
When swapping $r_i$ and $b_\ell$, we only apply local moves of the form $e_j f_j \leftrightarrow f_j e_j$.  
Now suppose $j=k$, which only happens when $k\leq \min(n-i, n-\ell)$.  In this case, by Proposition~\ref{prop:yd-minors} the minor goes from $([i,i+k-1], [\ell, \ell+k-1])$ to $([i+1, i+k], [\ell+1, \ell+k])$. 
However, the latter is a solid $k$-minor not present in the chamber minors of $Y_m$'s double wiring diagram, since it lies on the same diagonal as the former and has the same order. 
Thus it is in the test seed given by $Y_m$. By swapping this minor into the extended cluster in place of the original one, we get a new test seed. So, this mutation is a bridge. 
For $j>k$, the move $e_jf_j\leftrightarrow f_je_j$ corresponds to a mutation at a dead vertex.  
This mutation does not change the $k$-seed or $k$-positivity test at all, it only changes the way it is embedded in the total positivity cluster algebra.
Note that the size of $n^2$ is preserved since we never change the number of chambers in any track of the diagram, and the number of test variables is constant (since we only ever get new ones by swapping).
\end{proof}
\begin{definition}
We will call the paths described in Theorem~\ref{thm:young-dwd-correspondence} \emph{fundamental paths}.
\end{definition}

From the proof of this theorem, we can also easily prove the following fact:
\begin{corollary}
Each sub-cluster algebra found along a fundamental path has rank $(n-1)^2-(n-k)^2$.
\end{corollary}
\begin{proof}
The rank of the subcluster algebra is the number of mutable vertices in its quivers. The initial full $k$-quiver for the lexicographically minimal diagram has $(n-k)^2$ dead vertices and $2n-1$ frozen vertices (the $2(n-k)+1$ in track $k$ as well as the $2k-2$ unbounded chambers below track $k$), leaving $(n-1)^2-(n-k)^2$ mutable vertices. As discussed, no mutation at a dead vertex affects any of the mutable vertices.  A mutation at a chamber in track $k$, which occurs when jumping between sub-algebras, never adds edges between mutable and dead vertices since arrows only occur between chambers in adjacent tracks. Such mutations also always keep the frozen vertex adjacent to a dead one, which can be confirmed using conditions 2 and 3 from the definition of the quiver corresponding to a double wiring diagram. 
Therefore the number of mutable vertices is the same for every other quiver on the path.
\end{proof}

A path may travel through a number of components.  Any initial subsequence of a sequence as  described in Theorem~\ref{thm:young-dwd-correspondence} gives a $k$-positivity test preserving way to travel to some component that lies along the path corresponding to that sequence.  Notice from the proof of the theorem that we perform a bridge exactly when the groupings we are swapping both have a crossing in track $k$.
This means that every time the added box between $Y_{m}$ and $Y_{m+1}$ is placed inside the upper left $(n-k)\times (n-k)$ square of the Young diagram, a bridge occurs, and any time the added box is elsewhere, no bridge occurs, as all mutations specified by this box are at mutable vertices below track $k$. Also note that different Young diagrams within the $(n-k)\times (n-k)$ square give different components, since by Proposition~\ref{prop:yd-minors} the $k$-minors present in the diagram are distinct. Thus there is a 1-1 correspondence between Young diagrams contained in an $(n-k)\times (n-k)$ square and the components found along fundamental paths.

\begin{remark}
When $k=2$, the cluster algebras found along the fundamental paths are among those researched by Chmutov, Jiradilok, and Stevens called \emph{double rim hook cluster algebras}~\cite{CJS}.  Double rim hook cluster algebras are indexed by a sequence of north and east steps.  If a double rim hook cluster algebra is indexed by a sequence that contains $n-2$ north steps and $n-2$ east steps, then the double rim hook can be embedded in an $n\times n$ grid.  Removing the first row and column, the double rim hook cuts out a Young diagram in the upper right of the grid that fits in an $(n-2)\times (n-2)$ square.  The cluster algebra obtained from this double rim hook is the same as the cluster algebra along the fundamental path indexed by this Young diagram.  In their paper, Chmutov, Jiradilok, and Stevens describe all cluster variables for these cluster algebras.
\end{remark}

\section*{Acknowledgements}
This research was carried out as part of the 2017 REU program at the School of Mathematics at University of Minnesota, Twin Cities. The authors are grateful for the support of NSF RTG grant DMS-1148634. The authors would like to thank Pavlo Pylyavskyy and Victor Reiner for their guidance and encouragement. They are also grateful to Elizabeth Kelley for her many helpful comments on the manuscript and to Neeraja Kulkarni, Joe Suk, and Ewin Tang for the useful discussions on $k$-positivity.

\appendix
\section{Code}

All code used can be found at \href{https://github.com/ewin-t/k-nonnegativity}{https://github.com/ewin-t/k-nonnegativity}. In particular, we have code for generating the exchange graphs of the sub-cluster algebras for $k\leq 2$ or $n\leq 3$.


\printbibliography

\end{document}